\documentclass{amsart}
\usepackage{latexsym}
\usepackage{amstext}
\usepackage{amsmath}
\usepackage{amssymb}
\usepackage{amsthm}
\usepackage{enumitem}
\usepackage{mathrsfs} 
\usepackage{url}
\usepackage{comment}
\usepackage[utf8]{inputenc}
\usepackage[T1]{fontenc}
\usepackage[colorlinks={true},linkcolor={red},citecolor={blue}, urlcolor={blue}]{hyperref}
\usepackage{stmaryrd}              
\usepackage{verbatim}
\usepackage{mathtools}
\usepackage{mdframed}
\usepackage[margin=3 cm]{geometry}

\theoremstyle{plain}
\newtheorem{theorem}{Theorem}[section]

\newtheorem{corollary}[theorem]{Corollary}
\newtheorem{proposition}[theorem]{Proposition}
\newtheorem{lemma}[theorem]{Lemma}

\theoremstyle{definition}
\newtheorem{definition}[theorem]{Definition}

\newtheorem{remark}[theorem]{Remark}

\newtheorem{example}[theorem]{Example}
\newtheorem{conjecture}[theorem]{Conjecture}

\def \kbar {\overline{k}}
\def \O {\mathcal{O}}
\def \II {\mathcal{I}}

\DeclareMathOperator{\Supp}{Supp}

\DeclareMathOperator{\codim}{codim}
\DeclareMathOperator{\Proj}{Proj}

\allowdisplaybreaks
\title[Schmidt's subspace theorem for algebraic points]{A generalized Schmidt's subspace theorem for algebraic points of bounded degree}
\author{GuanHeng Zhao}

\begin{document}
	
\begin{abstract}
Schmidt's subspace theorem in terms of Seshadri constants for closed subschemes in subgeneral position has already been richly developed. Building upon this framework, we establish a Schmidt subspace theorem for algebraic points of bounded degree with respect to numerically equivalent ample divisors, attaining the optimal factor known to date and deriving a generalized weighted version. Also, we settle a long-standing conjecture in this field. 
\end{abstract}
	

\maketitle
\tableofcontents
\section{Introduction}
	
\subsection{Schmidt's subspace theorem for closed subschemes}
	
Schmidt's subspace theorem is a fundamental component in the realm of higher-dimensional Diophantine approximation. In its classical form, it asserts that for a number field $k$, the set of points in $\mathbb{P}^n(k)$ satisfying a certain inequality involving linear forms is contained in a finite union of proper linear subspaces. Over the years, it has been extended by numerous researchers in several directions. Notable contributions have been made by Evertse and Ferretti \cite{ef_festschrift}, Corvaja and Zannier \cite{CZ}, and Ru and Vojta \cite{RV}. 
	
A recent version of this theorem, which involves closed subschemes in general position and their Seshadri constants relative to a given ample divisor $A$, has been established.
This version builds upon the work of Ru and Wang \cite{RW22} and Vojta \cite{vojta_subscheme}, who have explored the use of beta constants as an alternative to Seshadri constants. 
This particular formulation was proven in the recent paper \cite{HL21}.
	
\begin{theorem}[{\cite[Theorem 1.3]{HL21}}]\label{1.1}
Let $X$ be a projective variety of dimension $n$ defined over a number field $k$. Let $S$ be a finite set of places of $k$.  For each $v\in S$, let $Y_{0,v},\ldots, Y_{n,v}$ be closed subschemes of $X$ in general position, defined over $k$. Let $A$ be an ample Cartier divisor on $X$ and $\epsilon>0$. Then there exists a proper Zariski-closed subset $Z\subset X$ such that for all points $P\in X(k)\setminus Z$,
\begin{equation*}
\sum_{v\in S}\sum_{i=0}^n \epsilon_{Y_{i,v}}(A) \lambda_{Y_{i,v},v}(P)< (n+1+\epsilon)h_A(P).
\end{equation*}
\end{theorem}

We first review the definition of the Seshadri constant of a closed subscheme $Y$ relative to a nef Cartier divisor $A$, which measures the local positivity of $A$ along $Y$.
	
Let $X$ be a projective variety and $Y$ be a closed subscheme of $X,$ corresponding to a coherent sheaf of ideals $\mathcal{I}_{Y}.$ Consider the sheaf of graded algebras $\mathcal{S} =
\oplus_{d \geq 0}\mathcal{I}_{Y}^{d}$, where $\mathcal{I}_{Y}^{d}$ is the $d$-th power of $\mathcal{I}_{Y},$ and $\mathcal{I}_{Y}^{0} = \mathcal{O}_{X}$ by convention. Then $\tilde{X} := \Proj \mathcal{S}$ is called the blowing-up of $X$ along $Y.$
	
Let $\pi: \tilde{X} \to X$ be the blowing-up of $X$ along $Y$. Note that the inverse image ideal sheaf $\tilde{\mathcal{I}}_Y = \pi^{-1} \mathcal{I}_Y \cdot \mathcal{O}_{\tilde{X}}$ is an invertible sheaf on $\tilde{X}$.

Then we have the following definition:
\begin{definition}\label{Sesh}
Let $Y$ be a closed subscheme of a projective variety $X$ and let $\pi:\tilde{X}\to X$ be the blowing-up of $X$ along $Y$. Let $A$ be a nef Cartier divisor on $X$. The {\it Seshadri constant} $\epsilon_Y(A)$ of $Y$ with respect to $A$ is the real number defined by
\begin{align*}
\epsilon_Y(A)=\sup\{\gamma\in {\mathbb{Q}}^{\geq 0}\mid \pi^*A-\gamma E\text{ is $\mathbb{Q}$-nef}\},
\end{align*}
where $E$ is an effective Cartier divisor on $\tilde X$ whose associated invertible sheaf is the dual of $\pi^{-1}\II_Y \cdot \O_{\tilde X}$.
		
In particular, when $Y=D$ is an effective Cartier divisor, we have
\begin{align}\label{S D}
\epsilon_{D}(A)=\sup\{\gamma\in {\mathbb{Q}}^{\geq 0}\mid A-\gamma D\text{ is $\mathbb{Q}$-nef}\}.
\end{align}
\end{definition}

Furthermore, we present Remark \ref{lem:seshadri_numerical} below, which will be useful in subsequent proofs and results, which incorporates the relationship between Seshadri constants and the coefficients in the numerically equivalent condition.
\begin{remark}\label{lem:seshadri_numerical}
Let $D$ be an effective Cartier divisor on a projective variety $X$.
Suppose $D\equiv dA$ for some positive integer $d$, where $A$ is an ample Cartier divisor. Then the Seshadri constant of $D$ with respect to $A$ is \[\epsilon_D(A)=\frac{1}{d}.\]
Indeed, by Definition \ref{Sesh}, for an effective Cartier divisor its Seshadri constant is given by (\ref{S D}).
Since the nef property is numerical and $D\equiv dA$, we have
$$A-\gamma D\equiv (1-\gamma d)A.$$ As $A$ is ample, 
$(1-\gamma d)A$ is nef if and only if $1-\gamma d\geq 0$. 
Taking the supremum over all such $\gamma$ yields 
$\epsilon_D(A)=1/d$.
\end{remark}
	
More recently, in \cite{Quang_Pacific}, Quang introduced the concept of a distributive constant for a family of divisors, and later in \cite{Quang22}, Quang also expanded this notable notion to encompass closed subschemes.
	
\begin{definition}\label{1.4}
For a family of closed subschemes $\mathcal{Y}=\{Y_1,\ldots, Y_q\}$ of a projective variety $X$, the distributive constant is given by
\begin{align*}
\delta_{\mathcal{Y}}=\max_{J\subset \{1,\ldots, q\}}\max\left\{1, \frac{\sharp J}{\codim \bigcap_{j\in J}\Supp Y_j}\right\},
\end{align*}
where we set $\codim\varnothing=\infty$.
\end{definition}
This concept contains all notions of subgeneral position:
\begin{definition}\label{1.5} Let $X$ be a projective variety of dimension $n$ defined over a number field $k$ and let $Y_{1}, \ldots, Y_{q}$ be $q$ closed subschemes of $X.$ Suppose $m$ and $\kappa$ are positive integers satisfying $m\geq\dim X\geq \kappa.$\par
		
(a). Closed subschemes $\{Y_{1}, \ldots, Y_{q}\}$  are said to be in general position in $X$ if for any subset $I\subset\{1, \ldots, q\}$ with $\sharp I\leq \dim X+1,$
$$\codim\left(\bigcap_{i\in I} Y_{i}\cap X\right)\geq\sharp I.$$\par
		
(b). Closed subschemes $\{Y_{1}, \ldots, Y_{q}\}$ are said to be in $m$-subgeneral position in $X$ if for any subset $I\subset\{1, \ldots, q\}$ with $\sharp I\leq m+1,$
$$\dim\left(\bigcap_{i\in I}Y_{i}\cap X\right)\leq m-\sharp I.$$\par
		
(c). Closed subschemes $\{Y_{1}, \ldots, Y_{q}\}$ are in $m$-subgeneral position in $X$ with index $\kappa$ if  $\{Y_{1}, \ldots, Y_{q}\}$  are in $m$-subgeneral position and for any subset $I\subset\{1, \ldots, q\}$ with $\sharp I\leq \kappa,$
$$\codim\left(\bigcap_{i\in I} Y_{i}\cap X\right)\geq\sharp I.$$
\end{definition}
It also yields the following important generalization.
\begin{theorem}[Quang \cite{Quang22}]
\label{1.6}
Let $X$ be a projective variety of dimension $n$ defined over a number field $k$, and let $S$ be a finite set of places of $k$. For each $v\in S$, let $Y_{1,v},\ldots, Y_{q,v}$ be closed subschemes of $X$, defined over $k$. Let $\Delta=\max_{v\in S}\delta_{\mathcal{Y}_v}$ be the maximum of the distributive constants of $\mathcal{Y}_v=\{Y_{1,v},\ldots, Y_{q,v}\}$. Let $A$ be an ample Cartier divisor on $X$, and $\epsilon>0$. Then there exists a proper Zariski closed set $Z$ of $X$ such that
\begin{align*}
\sum_{v\in S}\sum_{i=1}^q \epsilon_{Y_{i,v}}(A) \lambda_{Y_{i,v},v}(P)<  \left(\Delta(n+1)+\epsilon\right)h_A(P)
\end{align*}
for all points $P\in X(k)\setminus Z$.
\end{theorem}
	
As discussed in \cite{Quang_Pacific}, 
\begin{align*}
\Delta\leq m-n+1
\end{align*}
if the $Y_{i,v}$ are taken to be closed subschemes in $m$-subgeneral position.
\begin{align*}
\Delta\leq\frac{m-n}{\kappa}+1
\end{align*}
if the $Y_{i,v}$ are taken to be closed subschemes in $m$-subgeneral position with index $\kappa$.
	
Therefore, Theorem \ref{1.6} generalizes and encompasses all the prior results concerning closed subschemes in $m$-subgeneral position and their variations.

In addition, in the latest paper \cite{HL25}, the authors have extended Theorem \ref{1.6} to the following Schmidt-Nochka-type theorem for weighted sums involving arbitrary closed subschemes. 
	
\begin{theorem}[{\cite[Theorem 1.2]{HL25}}]
\label{1.7}
Let $X$ be a projective variety of dimension $n$  defined over a number field $k$, and let $S$ be a finite set of places of $k$. For each $v\in S$, let $Y_{1,v},\ldots, Y_{q,v}$ be closed subschemes of $X$ defined over $k$, and let $c_{1,v},\ldots, c_{q,v}$ be nonnegative real numbers.  For a closed subset $W\subset X$ and $v\in S$, let
\begin{align*}
\alpha_v(W)=\sum_{\substack{i\\ W\subset \Supp Y_{i,v}}}c_{i,v}.
\end{align*}
Let $A$ be an ample Cartier divisor on $X$, and $\epsilon>0$.  Then there exists a proper Zariski closed subset $Z$ of $X$ such that
\begin{align}\label{HA}
\sum_{v\in S}\sum_{i=1}^q c_{i,v}\epsilon_{Y_{i,v}}(A) \lambda_{Y_{i,v},v}(P)<  \left((n+1)\max_{\substack{v\in S\\ \varnothing\subsetneq W\subsetneq X}}\frac{\alpha_v(W)}{\codim W}+\epsilon\right)h_A(P)
\end{align}
for all points $P\in X(k)\setminus Z$.
\end{theorem}
	
Previous work on Schmidt's theorem relied primarily on the filtration technique. Theorem \ref{1.7} makes substantial use of weighted inequalities and Seshadri constants of intersections. This provides a new approach for addressing a variety of situations beyond the reach of the traditional filtration technique. Even when $c_{i,v}=1$ for all $i$ and $v \in S$, this theorem remains highly significant, mainly because of the form of the bound on the right-hand side.
	
In comparison with Theorem \ref{1.6}, under the assumption $c_{i,v}=1$ for all $i$ and $v$, we assert that
\begin{align*}
\delta_{\mathcal{Y}_v}=\max\left\{1, \max_{\varnothing\subsetneq W\subsetneq X}\frac{\alpha_v(W)}{\codim W}\right\}.
\end{align*}
This conclusion stems directly from the definitions and a straightforward observation, which we will elaborate on in  $\S$\ref{P3}, following the detailed explanation in \cite{CY}. When taking the maximum over $W$, it suffices to consider the cases where $W$ is the intersection of a subset of closed subschemes of $\mathcal{Y}_v$.
	
Overall, as stated in \cite{HL25}, Theorem \ref{1.7} can be regarded as a weighted extension of Theorem \ref{1.6}. Specifically, in certain cases, the formula (\ref{HA}) improves upon the coefficient of the height compared to Quang's theorem, primarily due to the presence of $1$ in the maximum in Quang's definition of $\delta_{\mathcal{Y}}$. Meanwhile, it significantly reduces the requirements on the positional condition. This furnishes our subsequent constructions with a decisive theoretical advantage.
	
\subsection{Theorems of algebraic points of bounded degree}
\label{R1.2}
While the above extensions concern rational points, a natural question is to consider algebraic points of bounded degree. For a point $P \in X(\kbar)$, let $[k(P):k]$ denote its degree over $k$. One asks whether the subspace theorem remains valid for points with $[k(P):k] \leq \delta$, where $\delta$ is a fixed positive integer.
In recent developments in this direction, our efforts have been consistently directed toward the identification of an optimal factor to bound the right-hand side of the inequality, especially for effective Cartier divisors and, as a special case, hyperplanes. Previous work studying Schmidt's theorem for algebraic points of bounded degree proved the following theorem by using the filtration technique. Henceforth, we write $\bigcup_{i,v}\Supp D_{i,v}=\bigcup_{v\in S}\bigcup_{i \in I}\Supp D_{i,v}$, where $D_{i,v}$ are effective Cartier divisors and $I$ is an index set.
\begin{theorem}[\cite{levin_duke}, Theorem 1.7]
\label{1.8}
Let $X$ be a projective subvariety of $\mathbb{P}^N$ with dimension $n\geq 1$ defined over a number field $k$.  Let $S$ be a finite set of places of $k$. For $v\in S$, let $H_{0,v},\ldots,H_{n,v}\subset \mathbb{P}^N$ be hypersurfaces over $k$ in n-subgeneral position (general position). Let $\delta\geq 1$ be an integer and $\epsilon>0$. Then the inequality
\begin{equation}\label{ineq1.8}
\sum_{v\in S}\sum_{\substack{w\in M_{k(P)}\\w\mid v}}\sum_{i=0}^n \frac{\lambda_{H_{i,v},w}(P)}{\deg H_{i,v}}< \left((\delta n)^2\left(\frac{\delta n-1}{2\delta n-3}\right)+\epsilon\right)h(P)
\end{equation}
holds for all but finitely many points $P\in X(\kbar)\setminus \bigcup_{i,v}H_{i,v}$ satisfying $[k(P):k]\leq \delta$.
\end{theorem}
Take $X = \mathbb{P}^n = \mathbb{P}^N$ and let $H_{i,v} \subset X$ be hyperplanes. Schlickewei (\cite{Schlickewei-03}) conjectured that there is a constant $c(\delta, n)$, depending only on $\delta$ and $n$, and a finite union of hyperplanes $Z \subset \mathbb{P}^n$ such that the left-hand side of inequality \eqref{ineq1.8} is bounded by $(c(\delta, n) + \epsilon)h(P)$ for all points $P \in \mathbb{P}^n(\kbar) \setminus Z$ satisfying $[k(P):k] \leq \delta$. Consequently, Theorem \ref{1.8} confirms Schlickewei's conjecture as a special case, but this does not give a sharp factor for sufficiently large $\delta$ and $n$.
	
Now, we consider more complex conditions. First, we explore some necessary parameters to better illustrate our statements in the following parts.
Consider a projective variety $X$ of dimension $n$ defined over a number field $k$, with a finite set of places $S$ of $k$. Let $m\geq n$ be a positive integer. For each place $v \in S$, we have $m+1$ effective Cartier divisors $D_{0,v}, D_{1,v}, \ldots, D_{m,v}$ on $X$ that are in $m$-subgeneral position. We define a divisor $D$ as the sum of all these divisors over all places in $S$,
$$D = \sum_{v \in S} \sum_{i=0}^m D_{i,v}.$$
Additionally, we are given an ample Cartier divisor $A$ on $X$ and positive integers $d_{i,v}$ such that each divisor $D_{i,v}$ is numerically equivalent to $d_{i,v}A$, that is, $D_{i,v} \equiv d_{i,v}A$ for all $i$ and $v \in S$. For $c\in\mathbb{R}$ and $\delta\in \mathbb{Z}^{+}$, define the geometric exceptional set $Z=Z\left(c,\delta, m, n,k, S, X, A, \{D_{i,v}\}\right)$ to be the smallest Zariski-closed subset of $X$ such that the set of points
\begin{equation}
\label{exc1}
\left\{P\in X(\kbar)\setminus \Supp D\mid 	\sum_{v\in S}\sum_{\substack{w\in M_{k(P)}\\w\mid v}}\sum_{i=0}^m \frac{\lambda_{D_{i,v},w}(P)}{d_{i,v}}>c h_A(P)\right\}\setminus Z
\end{equation}
satisfying $[k(P):k]\leq \delta$ is finite for all choices of the height functions.

Now define
\begin{multline*}
C(\delta, l,m,n)=\sup \{c\mid \dim Z\left(c,\delta, m, n,k, S, X, A, \{D_{i,v}\}\right)\geq l \\
\text{ for some $k$, $S$, $X$, $A$, $\{D_{i,v}\}$ as above}\}.
\end{multline*}

Furthermore, it would be beneficial for us to define and examine specific instances concerning rational points. These cases contain the initial factor without $\delta$. We define
\begin{align*}
A(m,n)=C(1,n,m,n), B(l,m,n)=C(1,l,m,n).
\end{align*}
Combined with the inequality (\ref{CBA}), which we will prove and explain in great detail (i.e., the proof of the relatively complex side of the index in \cite[Theorem 4.1]{levin_duke}), if $\delta m\geq 2$, there is a more general theorem, which admits a factor of
$$
C(\delta, 1,m,n)\leq \max_{1\leq j\leq \delta n}A(\delta m,j)\leq \max_{1\leq j\leq \delta n}\frac{\delta m(\delta m-1)(j+1)}{\delta m+j-2}=\frac{\delta m(\delta m-1)(\delta n+1)}{\delta m+\delta n-2}.
$$
Explicitly, this provides a version of Schmidt’s theorem for algebraic points of bounded degree.
\begin{theorem}[\cite{levin_duke}, Theorem 5.2]
\label{1.9}
Let $X$ be a projective variety of dimension $n$ defined over a number field $k$. Let $S$ be a finite set of places of $k$. Let $D_{0,v},\ldots, D_{m,v}$ be effective ample Cartier divisors on $X$ in $m$-subgeneral position for each $v\in S$. Suppose that there exists a Cartier divisor $A$ on $X$, defined over $k$, and positive integers $d_{i,v}$ such that $D_{i,v}\equiv d_{i,v}A$ for all $i$ and for all $v\in S$.  Let $\delta$ be a positive integer with $\delta m\geq 2$, and $\epsilon>0$. Then the inequality
\begin{equation*}
\sum_{v\in S}\sum_{\substack{w\in M_{k(P)}\\w\mid v}}\sum_{i=0}^m \frac{\lambda_{D_{i,v},w}(P)}{d_{i,v}}< \left(\frac{\delta m(\delta m-1)(\delta n+1)}{\delta m+\delta n-2}+\epsilon\right)h_A(P)+O(1)
\end{equation*}
holds for all points $P\in X(\kbar)\setminus \bigcup_{i,v} \Supp D_{i,v}$ satisfying $[k(P):k]\leq \delta$, where the $O(1)$-term depends on $\epsilon$.
\end{theorem}
If $m=n\geq 2$, then using $A(\delta n,\delta n)\leq \delta n+1$ gives the following bound in Theorem \ref{1.8}. We have
\begin{equation*}
	C(\delta, 1,n,n)\leq \max_{1\leq j\leq \delta n}A(\delta n,j)\leq (\delta n)^2\frac{\delta n-1}{2\delta n-3},
\end{equation*}
which is essentially quadratic in $\delta n$ for sufficiently large $\delta$ and $n$.

Since the proof of Theorem \ref{1.9} relies on the traditional filtration technique, this approach still fails to yield an optimal factor for algebraic points of bounded degree. Specifically, when $\delta=1$, this theorem reduces to the classical Schmidt subspace theorem for effective Cartier divisors in $m$-subgeneral position with the factor 
\begin{align}\label{L's fac}
\frac{m(m-1)(n+1)}{m+n-2}.
\end{align}
Direct computation shows $\frac{m(m-1)(n+1)}{m+n-2}-(m-n+1)(n+1)=\frac{(n-2)(n^2-1)}{m+n-2} \geq 0.$ 

Hence (\ref{L's fac}) is weaker than the bound $\Delta(n+1)$ in Theorem \ref{1.6}, thus also weaker than the factor in Theorem \ref{1.7}. Once we prove that the above two bounds also apply to algebraic points of bounded degree, we obtain the optimized new version, which is the main thread of the whole discussion.
	
\subsection{Main results}
Henceforth, except in Section \ref{sp}, we primarily consider the case of $m$-subgeneral position with the weights $c_{i,v}=1$ for all $i$ and $v \in S$. This case serves as our main object of study. 

Instead of directly using the filtration technique, and motivated by the latest work on the Schmidt subspace theorem and in alignment with the symmetric power method, we provide a new approach to solving the problem of subspace theorem for algebraic points of bounded degree. This yields the following Theorem \ref{main}, which is one of our main results. It can be regarded as the sharpest known generalization concerning algebraic points of bounded degree in the context of $m$-subgeneral position, and it concerns the geometric and arithmetic properties of Weil functions. For hyperplanes in general position, an alternative formulation of this theorem, given in (\ref{4}), yields an optimal factor $2\delta n$ and thus settles a conjecture (see Section \ref{OC} for details). This value also strictly improves the previous bound (\ref{ineq1.8}), which was approximately $(\delta n)^2$.

\subsubsection{The factor $a_\delta$}\label{a1}\quad\par
In the sequel, we let $$\ell(m,n)=(n+1)\max_{\substack{v\in S\\ \varnothing\subsetneq W\subsetneq X}}\frac{\alpha_v(W)}{\codim W}.$$
For a positive integer $\delta$, define 
\[a_\delta=\max_{1\leq j\leq \delta n}\ell(\delta m,j),\]
here, $\ell(\delta m,j)$ denotes the same expression as above with the pair $(m,n)$ formally replaced by $(\delta m,j)$.
\begin{theorem}
\label{main}
Let $X$ be a projective variety of dimension $n$ defined over a number field $k$, and let $S$ be a finite set of places of $k$. For each $v\in S$, let $D_{1,v},\ldots, D_{m+1,v}$ be effective ample Cartier divisors on $X$, defined over $k$ and in $m$-subgeneral position. Let $d_{i,v}$ be positive integers such that $D_{i,v}\equiv d_{i,v}A$ for all $i$ and for all $v\in S$, where $A$ is an ample Cartier divisor defined on $X$. For a closed subset $W\subset X$ and $v\in S$, let
\begin{align*}
\alpha_v(W)=\sharp\left\{i \mid W\subset \Supp D_{i,v} \right\}.
\end{align*}
Let $\delta$ be a positive integer, and let $a_\delta$ be  as defined above. Let $\epsilon>0$. Then the inequality
\begin{align*}
\sum_{v\in S}\sum_{\substack{w\in M_{k(P)}\\w\mid v}}\sum_{i=1}^{m+1} \frac{\lambda_{D_{i,v},w}(P)}{d_{i,v}}<  \left(a_\delta+\epsilon\right)h_A(P)+O(1)
\end{align*}
holds for all points $P\in X(\kbar)\setminus \bigcup_{i,v} \Supp D_{i,v}$ satisfying $[k(P):k]\leq \delta$, where the $O(1)$-term depends on $\epsilon$.
\end{theorem}

\begin{remark} [Implication of $a_\delta$]
\normalfont
\label{1.11}
Actually, it is an expression that can be obtained by reducing Diophantine approximation problems for algebraic points to those for rational points on symmetric powers of $X$.
An overview of the steps is as follows:
\begin{enumerate}[label=(\Roman*)., leftmargin=2.5em, nosep]
	\item The rational point case ($\delta=1$). 
	
	For rational points, we prove the initial factor $\ell(m,n)$ with respect to numerically equivalent ample divisors. This confirms that $\ell(m,n) \leq \Delta(n+1) \leq (m-n+1)(n+1) \leq \frac{m(m-1)(n+1)}{m+n-2}.$
	
	\item Algebraic points of bounded degree $[k(P):k]\leq \delta$. 
	
	We employ the symmetric power construction as a central tool (see Section \ref{section_main_thm_proof}, Part II) to handle this problem: Let $X^\delta$ denote the $\delta$-th power of $X$. We replace the base variety $X$ by its $\delta$-th symmetric power $X^{(\delta)}$, replace each divisor $D_{i,v}$ by $D_{i,v}^*=\phi_*\pi_1^*(D_{i,v})$, and replace an ample divisor $A$ by $A^* = \phi_* \pi_1^* (A)$, where $\phi:X^\delta \to X^{(\delta)}$ is a quotient map and $\pi_i: X^\delta \to X$ denote the $i$-th projection maps. Then the new dimension is $\delta n$ and new subgeneral position parameter is $\delta m$. Overall, we can thus provide a complete proof of the inequality (\ref{CBA}).
	
	\item For any function $\mathcal{L}(m,n)$ satisfying the given conditions and defined in the context of rational points, we denote its degree-$\delta$ extension by $$\max_{1\leq j\leq \delta n}\mathcal{L}(\delta m,j).$$
	This operation consists of\par
	\begin{itemize}
	\item Verifying the new parameters $(\delta m, j)$ after symmetric power construction;
	
	\item Maximizing over all possible exceptional component dimensions $j \in \left\{1,\ldots,\delta n\right\}$ to ensure uniformity.
    \end{itemize}
	Also, the inequality mentioned above confirms the existence of the factor. We have
	\begin{align*}
		&C(\delta, 1,m,n)\leq \max_{1\leq j\leq \delta n}A(\delta m,j)\leq \max_{1\leq j\leq \delta n}\mathcal{L}(\delta m,j).
	\end{align*}
	
	\item Applying this operation to $\ell(m,n)$, we obtain the final expression
	$$a_\delta=\max_{1\leq j\leq \delta n}\ell(\delta m,j).$$
	
\end{enumerate}
According to \cite[Section 4 $\&$ Section 5]{levin_duke}, this is precisely the result we need, as it allows us to convert the rational point factor into the algebraic point factor. Moreover, since the function $\ell(m,n)$ gives a sharper bound than any previous result, $a_\delta$ inherits this property directly, thereby guaranteeing that our result is the sharpest to date. That is
$$a_\delta \leq \frac{\delta m(\delta m-1)(\delta n+1)}{\delta m+\delta n-2}.$$

\end{remark}

We now list some significant results for hypersurfaces and hyperplanes. Also, in the rest of this section, we obtain some slightly more general formulations for algebraic points of bounded degree (e.g., with index $\kappa$, see the next subsection \S\ref{a2}), which constitute a generalized extension.
	
Let $X \subset \mathbb{P}^N(k)$ be a projective variety of dimension $n$ defined over a number field $k$. If $D$ is a hypersurface in $\mathbb{P}^N(k)$ defined by a homogeneous polynomial $f\in k[x_0,\ldots,x_N]$ of degree $d$, then for $v\in M_k$, let $|f|_v$ denote the maximum of the absolute values of the coefficients of $f$ with respect to $v$. We define $\|f\|_v$ similarly. Then a Weil function for $D$ is given by
\begin{equation*}
	\lambda_{D,v}(P)=\log \frac{\|f\|_v\max_{i} \|x_i\|_v^d}{\|f(P)\|_v}, \quad 0\leq i \leq N.
\end{equation*}
\begin{corollary}\label{p}
Let $\mathcal{D}_v=\{D_{1,v},\ldots, D_{m+1,v}\}$ be hypersurfaces in $m$-subgeneral position. Then
\begin{align*}
\log\prod_{v\in S}\prod_{\substack{w\in M_{k(P)}\\w\mid v}}\prod_{i=1}^{m+1} \left(\frac{\|f_{i,v}\|_{w}\max_j \|x_{j}\|^{d_{i,v}}_{w}}{\|f_{i,v}(P)\|_{w}}\right)^{\frac{1}{d_{i,v}}}<\left(a_\delta+\epsilon\right)h_A(P)+O(1)
\end{align*}
holds for all points $P\in X(\kbar)\setminus \bigcup_{i,v} \Supp D_{i,v}$ satisfying $[k(P):k]\leq \delta$, where $A$ is a hyperplane.
\end{corollary}

Now assume that $\mathcal{D}=\{D_{1},\ldots, D_{q}\}$ is a set containing $q$ (not necessarily $m+1$ here, as confirmed by Proposition \ref{p4.2} and Remark \ref{4.3}) hyperplanes in $m$-subgeneral position and consider a setup as stated in Section \ref{OC}. Then the inequality
\begin{align}\label{4}
	\log\prod_{v\in S}\prod_{\substack{w\in M_{k(P)}\\w\mid v}}\prod_{i=1}^q\left(\frac{\|f_{i}\|_{w}\max_j \|x_{j}\|_{w}}{\|f_{i}(P)\|_{w}}\right)< \left(2\delta m+\epsilon\right)h(P)+O(1)
\end{align}
holds for all points $P\in X(\kbar)\setminus \Supp\sum_{i=1}^{q} D_i$, where the $O(1)$-term depends on $\epsilon$.

If the situation specializes to general position, so that $m=n \geq 2$, then we arrive at the following inequality 
$$C(\delta,1,n,n) \leq B(1,\delta n, \delta n) \leq \max_{1\leq j\leq \delta n}A(\delta n,j) \leq \max_{1\leq j\leq \delta n}\left(\frac{\delta n}{j}\right)(j+1)=2\delta n=a_\delta.$$

In 2000, Ru proposed the following conjecture about Wirsing’s version of Schmidt’s subspace theorem for linear forms (which can also be found in \cite[Conjecture 7.1, part (b)]{levin_duke}). This conjecture in number theory is motivated by the correspondence in Nevanlinna theory over complex geometry on the covering space in $\mathbb{C}^n$.

For each $v \in M_k$, we let $k_v$ denote the completion of $k$.
We set $\|x\|_v=|x|_v^{[k_v:\mathbb{Q}_v]/[k:\mathbb{Q}]}.$

\begin{conjecture}[{\cite[Conjecture 4.4 with $l=1$]{RVC}}]
\label{c}
Let $k$ be a number field, let $M_k$ be the set consisting of all the
non-equivalent places of $k$. Let $S \subset M_k$ be a finite set containing all the archimedean places. Let $\{L_{1},\ldots, L_{q}\}$ be linear forms with coefficients in $k$, in general position. Let $\delta$ be a positive integer such that $[k(\mathbf{x}):k]\leq \delta$. If $\epsilon > 0$, and if $C$ is a positive constant, then the set of points $\mathbf{x} \in [x_0:x_1:\cdots:x_n]$ with $x_j \in \kbar, 0 \leq j \leq n$, such that 
\begin{align*}
\log\prod_{v\in S}\prod_{\substack{w\in M_{k(\mathbf{x})}\\w\mid v}}\prod_{i=1}^q \left(\frac{\|L_{i}\|_{w}\max_j \|x_{j}\|_{w}}{\|L_{i}(\mathbf{x})\|_{w}}\right)\geq\left(2\delta n+\epsilon\right)h(\mathbf{x})+C
\end{align*}
is contained in a finite union of linear subspaces of dimension $0$.
\end{conjecture}
\begin{remark}
	In the setting of hyperplanes, i.e. hypersurfaces of degree $d=1$, the Weil function defined in Corollary \ref{p} reduces to the classical linear form height
	$$\lambda_{D,v}(P)=\log \frac{\|L\|_v\max_{i} \|x_i\|_v}{\|L(P)\|_v}, \quad 0\leq i \leq N.$$
	where $L$ is a linear form defining the hyperplane $D$.
\end{remark}
For $l=1$, this special case seems representative and is likely to be quite difficult already via the traditional filtration technique. Our new approach will overcome these obstacles by using the internal properties of Weil functions corresponding to the intersection of closed subschemes yielding an explicit and optimal factor.

For the algebraic variety $X=\mathbb{P}^n=\mathbb{P}^N$, in accordance with Example \ref{EX}, we will prove that $$C(\delta,1,n,n)=B(1,\delta n, \delta n)=2\delta n.$$ 
We recall the definition of 
\begin{multline*}
	C(\delta, 1,n,n)=\sup \{c\mid \dim Z\left(c,\delta, n, n,k, S, X, A, \{D_{i,v}\}\right)\geq 1 \\
	\text{ for some $k$, $S$, $X$, $A$, $\{D_{i,v}\}$ as in (\ref{exc1})}\}.
\end{multline*}
Since $C(\delta,1,n,n)=2\delta n$, for any constant $c>2\delta n$ the exceptional set $Z$ satisfies $\dim Z<1$, hence $Z$ is a finite set of points and in particular is contained in a finite union of linear subspaces of dimension $l-1=0$. As a consequence, we settle Conjecture \ref{c}.

\subsubsection{The factor $a_{\Delta,\delta}$}\quad\par
\label{a2}
Also, by Theorem \ref{1.7} and Remark \ref{1.11}, we slightly weaken the strength of the factor $a_\delta$. 
We replace $a_\delta$ by a more flexible and explicit factor $a_{\Delta,\delta}$ (with $a_\delta \leq a_{\Delta,\delta}$),
defined as
$$\max_{1\leq j\leq \delta n}\iota(\delta m,j),$$
where $\iota(m,n)=(m-n+1)(n+1)$, or $(\frac{m-n}{\kappa}+1)(n+1)$, associated with distributive constants $\delta_{\mathcal{D}_v}$. 

Although it sacrifices some sharpness, this formulation is geometrically more refined and allows a generalization to $m$-subgeneral position with index $\kappa$. Therefore, in the following discussion, we can equip this factor with the index $\kappa$. We now state a proposition involving this new expression.

\begin{proposition}\label{pp}
	Under the same assumptions as in Theorem \ref{main}, for each $v \in S$, 
	let $\mathcal{D}_v=\{D_{1,v},\ldots, D_{m+1,v}\}$ be a set of effective ample Cartier divisors (in $m$-subgeneral position or with index $\kappa$). Let $\epsilon>0$ and let $\delta$ be a positive integer with $\delta m > 2$. Then the inequality
	\begin{align*}
		\sum_{v\in S}\sum_{\substack{w\in M_{k(P)}\\w\mid v}}\sum_{i=1}^{m+1} \frac{\lambda_{D_{i,v},w}(P)}{d_{i,v}}<  \left(a_{\Delta,\delta}+\epsilon\right)h_A(P)+O(1)
	\end{align*}
	holds for all points $P\in X(\kbar)\setminus \bigcup_{i,v} \Supp D_{i,v}$ satisfying $[k(P):k]\leq \delta$, where the $O(1)$-term depends on $\epsilon$.
\end{proposition}

At this point, we have, for $m$-subgeneral position, $\iota(m,n)=(m-n+1)(n+1)$ and
$$a_{\Delta,\delta}=\max_{1\leq j\leq \delta n}\left(\delta m-j+1\right)(j+1).$$	
Then we conclude the following main results:

(i) For $m>2n$,
$$a_{\Delta,\delta}={\delta^2}n(m-n)+{\delta}m+1 .$$

(ii) For $n \leq m \leq 2n$, let $*$ be a positive real number and $[*]$ denote the rounding symbol, taking the nearest positive integer to $*$. It is easy to see that $[*]$ equals $*$ if this number is already an integer, and when $*$ is exactly the midpoint, two integers satisfy the condition.

Hence, for $j$ ($1 \leq j \leq \delta n$) defined as above,  the maximum occurs when $j=\left[\frac{\delta m}{2}\right]$. We denote this index by $i$, thus $$a_{\Delta,\delta}=-i^2+\delta m \cdot i+{\delta}m+1 .$$
Similarly, in higher dimensions and for larger $\delta$, this new bound $a_{\Delta, \delta}+\epsilon$ also significantly improves the expression $$\frac{\delta m(\delta m-1)(\delta n+1)}{\delta m+\delta n-2}+\epsilon$$
in Theorem \ref{1.9}.
	
\begin{remark}
It can be verified that the precise difference between the two results is $$\frac{(\delta n-2)(\delta^2 n^2-1)}{\delta m+\delta n-2},$$ in the setting of $m>2n$. For $n \leq m \leq 2n$, see the calculation in Section \ref{P3}.
\end{remark}
Furthermore, for $\mathcal{D}_v$ in $m$-subgeneral position with index $\kappa$, let $\iota(m,n)=(\frac{m-n}{\kappa}+1)(n+1)$. Then
$$a_{\Delta,\delta}=\max_{1\leq j\leq \delta n}\left(\frac{\delta m-j}{\kappa}+1\right)(j+1),$$ 
and we can conclude that\par
(i) For $m \geq \lceil 2n+\frac{1-\kappa}{\delta} \rceil$,
$$a_{\Delta,\delta}=\frac{n}{\kappa}\left({\delta^2}m-{\delta^2}n+\delta(\kappa-1)\right)+\frac{{\delta}m}{\kappa}+1 .$$\par 
(ii) For $n \leq m \leq  \lfloor 2n+\frac{1-\kappa}{\delta} \rfloor $, in the same way, we have $i=\left[\frac{\delta m+(\kappa-1)}{2}\right],$ thus
$$a_{\Delta,\delta}=-\frac{i^2}{\kappa}+\left(\frac{\delta m+(\kappa -1)}{\kappa}\right)i+\frac{{\delta}m}{\kappa}+1 .$$
In this generalized case, since the definition of $a_\delta$ does not directly involve $\kappa$, the factor $a_{\Delta, \delta}$ is more explicit and yields a new result that contains the notion of index $\kappa$.

The structure of this paper is as follows.
In Section \ref{sec_def_prelim}, we introduce some notions and recall some basic definitions for height functions. 
In Section \ref{section_main_thm_proof}, we divide the problem into three parts, along the lines of \cite{levin_duke}. We first provide a rigorous proof of Schmidt's subspace theorem for effective Cartier divisors based on the idea in \cite{HL25} and extend it to numerically equivalent ample divisors in two different ways, and then handle the case of algebraic points of bounded degree. In Section \ref{sp}, we present some supplementary propositions that are weakened variants derived from the proof of the main theorem. In Section \ref{OC}, we discuss a special case showing the optimality of our result under the conditions that $X=\mathbb{P}^n=\mathbb{P}^N$ and that the effective Cartier divisors $D_{i}$ are hyperplanes in general position.
	
\section{Preliminaries and height functions}\label{sec_def_prelim}

Throughout this paper, we denote by $k$ a number field, by $\kbar$ its algebraic closure, and by $\mathcal{O}_k$ its ring of integers. We define $M_k$ as a canonical set of places (or absolute values) of $k$, which includes:\par
(i) one place for each prime ideal $\mathfrak{p}$ of $\mathcal{O}_k$;\par
(ii) one place for each real embedding $\sigma: k \to \mathbb{R}$;\par
(iii) one place for each pair of conjugate embeddings $\sigma, \overline{\sigma}: k \to \mathbb{C}$.\par
Given a finite set $S$ of places of $k$ that contains all the archimedean places, we denote $\mathcal{O}_{k,S}$ as the ring of $S$-integers of $k$ and $\mathcal{O}_{k,S}^*$ as the group of $S$-units of $k$.
For a place $v$ of $k$ and a place $w$ of a field extension $L$ of $k$, we say that $w$ lies above $v$ (denoted $w|v$) if $w$ and $v$ induce the same topology on $k$. 

We normalize the absolute values as follows:
$|p|_v = \frac{1}{p}$ if $v$ corresponds to a prime ideal $\mathfrak{p}$ lying above a rational prime $p$,
$|x|_v = |\sigma(x)|$ if $v$ corresponds to an embedding $\sigma$. 
For each $v \in M_k$, we let $k_v$ denote the completion of $k$ with respect to $v$.
We set
\begin{equation*}
\|x\|_v=|x|_v^{[k_v:\mathbb{Q}_v]/[k:\mathbb{Q}]}.
\end{equation*}
A fundamental equation is the product formula
\begin{equation*}
\prod_{v\in M_k}\|x\|_v=1,
\end{equation*}
which holds for all $x\in k^*$.

In addition, we define the vector space $$L(D) = \{ f \in \kbar(X) \mid \text{div}(f) \geq -D \},$$ where $\text{div}(f)$ denotes the divisor associated with $f$. 

The dimension of this vector space over $\kbar$ is denoted by $l(D) = \dim_{\kbar} L(D)$. We use the symbol $\sim$ to indicate linear equivalence of divisors and $\equiv$ to denote numerical equivalence of divisors.	

Now we provide a concise overview of the essential properties of local height functions, primarily following the references \cite{silverman_87} and \cite{HL25}.
For more refined versions and more details about local or global height functions, see \cite[Section 2]{HL25}. We will only give a brief explanation here.
	
\subsection{Local height functions}\label{lhf}
	
Let $Y$ be a closed subscheme of a projective variety $X$, both defined over a number field $k$.  For any place $v$ of $k$, one can associate a {\it local height function} (or {\it Weil function}) $$\lambda_{Y,v}:X(k)\setminus Y\to \mathbb{R},$$ up to a bounded term $O(1)$. 

In particular, any $\lambda_{Y,v}$ is bounded from below for all $P\in X(k)\setminus Y$, and for local height functions ($\lambda_{Y,v})_{v \in M_k}$, they are also bounded from below by an $M_k$-constant (a sequence of real numbers $(c_v)_{v\in M_K}$ where $c_v = 0$ for all but finitely many $v$). Local height functions also satisfy the following properties. 

If $Y$ and $Z$ are two closed subschemes of $X$, defined over $k$, and $v$ is a place of $k$, then up to $O(1)$,
$$\lambda_{Y\cap Z,v}=\min\{ \lambda_{Y,v},\lambda_{Z,v}\}, \lambda_{Y+Z,v}=\lambda_{Y,v}+\lambda_{Z,v}, \lambda_{Y,v}\leq \lambda_{Z,v}, \text{ if }Y\subset Z.$$
It is worth highlighting that by setting $Y=Z$, we can obtain the following equality
$$\lambda_{2Y,v}= \lambda_{Y+Y,v}=\lambda_{Y,v}+\lambda_{Y,v}=2\lambda_{Y,v},$$
where $2Y$ denotes the closed subscheme corresponding to the ideal sheaf $\II_Y^2$. Furthermore, it is generally true that if $Y_1,\ldots,Y_r$ are closed subschemes of $X$ and $c_1,\ldots,c_r$ are positive integers, we have
\begin{align}
\lambda_{c_1Y_1\cap\ldots\cap c_rY_r,v}&=\min\{ c_1\lambda_{Y_1,v},\ldots,c_r\lambda_{Y_r,v}\},\label{cap}\\
\lambda_{c_1Y_1+\cdots+c_rY_r,v}&=c_1\lambda_{Y_1,v}+\cdots+c_r\lambda_{Y_r,v}.\label{plus}
\end{align}\par
If $\phi:W\to X$ is a morphism of projective varieties with $\phi(W) \not \subset Y$, then up to $O(1)$,
\begin{equation*}
\lambda_{Y,v}(\phi(P))=\lambda_{\phi^*Y,v}(P),
\quad \forall P\in W(k)\setminus \phi^*Y.
\end{equation*}
The operations $Y\cap Z$, $Y+Z$, $Y\subset Z$, and $\phi^*Y$ are defined in terms of the associated ideal sheaves. Specifically, it is a subtle point that if $Y$ is associated with the ideal sheaf $\mathcal{I}_Y$, then $\phi^* Y$ represents the closed subscheme corresponding to the inverse image ideal sheaf $\phi^{-1} \mathcal{I}_Y \cdot \mathcal{O}_W$. \par 
If $Y=D$ is an effective Cartier divisor, which we often identify with the associated closed subscheme according to \cite[Lemma 2.2]{silverman_87}, these height functions coincide with the usual height functions associated with divisors. 
	
\subsection{Global height functions} 
	
Let $Y$ and $X$ be as in subsection \ref{lhf}. We define a {\it global height function} $h_Y:X(k)\setminus Y\to \mathbb{R}$ associated with $Y$ as follows, 
$$h_Y(P) = \sum_{v\in M_k} \lambda_{Y,v} (P).$$
For $v\in M_k$, the definition of $\lambda_{Y,v}$ can be chosen such that the collection $(\lambda_{Y,v})_{v\in M_k}$ is well-defined up to an $M_k$-constant, thus $h_Y$ is well-defined up to $O(1)$. Global height functions possess properties similar to those described earlier for local height functions, with one notable exception: the property above for the local height function $\lambda_{Y\cap Z,v}=\min\{ \lambda_{Y,v},\lambda_{Z,v}\} $ is modified to $h_{Y \cap Z} \leq \min\{h_Y, h_Z\} + O(1).$

Particularly, there is an interpretation of arithmetic complexity behind global functions as follows. If a point $P$ is $v$-adically close to $Y$, then after a projective embedding using a very ample multiple of the divisor $A$, the coordinates of $P$ become $v$-adically complex. In \cite[Proposition 4.2]{silverman_87}, the author gave a concrete explanation of this statement. Additionally, this implies that for any arbitrary closed subscheme $Y$, $h_A$ dominates $h_Y$, as the domain of $h_A$ can be extended to contain all of $X(k)$ and there exists a constant $c$ such that for $P\in X(k)\setminus Y$, we have
\begin{align} 
h_Y(P)\leq c h_A(P). \label{ample_domination}
\end{align}
Also, given the properties of local height functions, the inequality \eqref{ample_domination} can be stated for each individual $\lambda_{Y,v}$: 
\begin{align} 
\lambda_{Y,v}(P)\leq c h_A(P). \label{lhf_ample_domination}
\end{align}
	
\section{Proof of Theorem \ref{main} and Proposition \ref{pp}}\label{section_main_thm_proof}
	
In this section, we divide our proof into the following implications.

\begin{itemize}
\item (a): Schmidt's theorem for numerically equivalent ample divisors
\end{itemize}
\begin{equation*}
\Downarrow \textbf{(i)}
\end{equation*}\par
\begin{itemize}
\item (b): (a) $+$ controlling the dimension of the exceptional set $Z$
\end{itemize}\par
Item \textbf{(i)} is achieved through repeated applications of (a) to exceptional sets.
\begin{equation*}
\Downarrow \textbf{(ii)}
\end{equation*}\par
\begin{itemize}
\item (c): (b) for algebraic points of bounded degree on $X$
\end{itemize}\par
Item \textbf{(ii)} is achieved by reducing Diophantine approximation problems for algebraic points to those for rational points on symmetric powers of $X$.
	
In brief, after having proved the theorem of type (a), we can derive the results stated in (b) and (c) immediately. As for the proof of type (a), our first method largely presents the strength and logical relation between several different conditions, even though it is more complicated.
	
\begin{remark}
According to the standard definition, the Weil function $\lambda_{D_{i,v}}(P)$ for an effective Cartier divisor  $D_{i,v}$ on $X$ is also determined by a continuous (fiber) metric $\|\cdot\|$ on the line sheaf $\mathcal{O}_X(D_{i,v})$, while $s_{i,v} \in H^0(X, \mathcal{O}_X(D_{i,v}))$ is its canonical section. For brevity we write $$\Omega_{D_{i,v}}(P)=\|s_{i,v}(P)\|^{-1}.$$
It is reasonable and meaningful to use the logarithmic representation of Weil functions (defined up to $O(1)$) to better elaborate on Corollary \ref{p}, Conjecture \ref{c}, and other results concerned with concrete metrics. We let  $$\epsilon_{D_{i,v}}(A)\lambda_{D_{i,v}}(P)=\log\Omega^{\epsilon_{D_{i,v}}(A)}_{D_{i,v}}(P)$$
throughout the discussion.
\end{remark}

In order to prove Theorem \ref{main}, we require the following inequalities from \cite{HL25} (The proof is analogous so we will not elaborate further here), which are a generalization of Chebyshev's inequality (see also work in \cite{MPF}). 
\begin{lemma}
\label{elemineq}
Let $\lambda_1\geq \lambda_2\geq\cdots \geq \lambda_n\geq 0$ be a family of Weil functions, write $\lambda_i=\log\Omega_i$ and let $b_1,\ldots, b_n,c_1,\ldots, c_n$ be nonnegative real numbers as sets of indicators. Let $i_0$ be the smallest index $i$ such that $c_{i}\not = 0$ (assume that there exists at least one such index). Then
\begin{align}
\label{generalChebyshev}
\log\prod_{i=1}^n\Omega^{b_i}_i\geq \log\prod_{i=1}^n\Omega^{\left(\min_{i_0\leq j\leq n} \frac{\sum_{i=1}^jb_i}{\sum_{i=1}^jc_i}\right)c_i}_i.
\end{align}
\end{lemma}
\begin{corollary}
\label{cor_elemineq}
Let $\lambda_1\geq \lambda_2\geq\cdots \geq \lambda_n\geq 0$ be a family of Weil functions, write $\lambda_i=\log\Omega_i$ and let $b_1,\ldots, b_n,c_1,\ldots, c_n$ be nonnegative real numbers as sets of indicators. Assume that $b_1\not = 0$. Then
\begin{align*}
\log\prod_{i=1}^n\Omega^{\left(\max_{1\leq j\leq n} \frac{\sum_{i=1}^jc_i}{\sum_{i=1}^jb_i}\right)b_i}_i\geq \log\prod_{i=1}^n\Omega^{c_i}_i.
\end{align*}
\end{corollary}

We now conduct the proof of Theorem \ref{main} and Proposition \ref{pp}.

\subsection{Part I}\label{P1}
Schmidt's theorem for numerically equivalent ample divisors.

In this part, let $m$ be a positive integer with $m \geq n$, and assume that $q \geq m+1$.

$\underline{{Approach}-\mathbb{I}}$: Let $v\in S$, and let $c_{i,v}$ $(1 \leq i \leq q)$ be nonnegative real numbers. Let $D_{i,v}$ be effective Cartier divisors.\par
\textbf{(i)} We normalize in a two-step process.

(a) Rationalization: First, for every $i=1\ldots, q$, let $\tau_{1,v},\ldots,\tau_{q,v}$ be real numbers in the interval $[0,1)$ so that 
$(1+\tau_{1,v})\epsilon_{D_{1,v}}(A),\ldots, (1+\tau_{q,v})\epsilon_{D_{q,v}}(A)$
are rational numbers. If $\epsilon_{D_{i,v}}(A)$ is already a rational number, then choose $\tau_{i,v}=0$.\par 
(b) Integer conversion: Fix a sufficiently divisible integer $c_v$ (independent of $i$) such that 
$$ c_v (1+\tau_{1,v})\epsilon_{D_{1,v}}(A),\ldots, c_v(1+\tau_{q,v})\epsilon_{D_{q,v}}(A)$$
are all integers.\par
Multiplying the effective Cartier divisors $D_{i,v}$, we set 
$$\tilde  D_{i,v} := c_v (1+\tau_{i,v})\epsilon_{D_{i,v}}(A)D_{i,v},$$
and obtain
\begin{equation}
\epsilon_{\tilde  D_{i,v}}(A)=\frac 1 {c_v(1+\tau_{i,v})} \label{modif_sesh}
\end{equation}
for $i=1,\ldots,q$, and therefore
\begin{equation*}
\frac 1 {c_v} \leq 2\epsilon_{\tilde  D_{i,v}}(A),
\end{equation*}
which implies
\begin{equation*}
\frac 1 {c_v} \leq 2 \min_{i}  \epsilon_{\tilde  D_{i,v}}(A).
\end{equation*}
For each Seshadri constant and for all $i,i'=1,\ldots,q$, we have 
\begin{align}
|\epsilon_{\tilde  D_{i,v}}(A)-\epsilon_{\tilde  D_{i',v}}(A)| 
\leq \max_k\tau_{k,v} \cdot 2  \min_{l}  \epsilon_{\tilde  D_{l,v}}(A). \label{sesh_bound}
\end{align} 
		
\textbf{(ii)} Intersections of the normalized divisors.

Set
$$\tilde D_{I,v}=\bigcap_{i\in I}\tilde  D_{i,v}.$$ 
In fact, $\tilde{D}_{I,v}$ are defined to be closed subschemes by using \cite[Lemma 2.2]{silverman_87}. We simply refer to them as divisors in the sequel. Fix a point $P$ outside of the support of the given divisors. Let $\{i_{1,v},\ldots, i_{q,v}\}=\{1,\ldots, q\}$ be such that:
\begin{align}
\log\Omega_{\tilde D_{i_{1,v},v},v}(P)\geq \cdots \geq \log\Omega_{\tilde D_{i_{q,v},v},v}(P)\geq 0, \label{tilde_chain_ineq}
\end{align}
and set $I_{j,v}=\{i_{1,v},\ldots, i_{j,v}\}$.  Let $l_v$ be the largest index such that
$\bigcap_{j=1}^{l_v} \tilde D_{i_{j,v},v}\neq \varnothing.$

Using property \eqref{cap} and the sequence of inequalities \eqref{tilde_chain_ineq}, we have
\begin{align*}
\log\Omega_{\tilde D_{i_{j,v}v},v}(P)=\min_{j'\leq j}\log\Omega_{\tilde D_{i_{j',v},v},v}(P)=\log\Omega_{\tilde D_{I_{j,v},v},v}(P).
\end{align*}
In particular, for $j>l_v$, $\log\boxdot_{\tilde D_{i_{j,v},v},v}(P)$ is bounded above by a constant independent of $P$, and
\begin{align*}
\log\prod_{v\in S}\prod_{i=1}^q \Omega^{c_{i,v}\epsilon_{\tilde D_{i,v}}(A)}_{\tilde D_{i,v},v}(P)\leq \log\prod_{v\in S}\prod_{j=1}^{l_v} \Omega^{c_{i_j,v}\epsilon_{\tilde D_{i_{j,v},v}}(A)}_{\tilde D_{I_{j,v},v},v}(P). 
\end{align*}
Let $b_{j,v}=\codim \tilde D_{I_{j,v},v}$ and set $b_{0,v}=0$.  Set $\epsilon_v(A)=\max_i \epsilon_{\tilde D_{i,v}}(A)$.  Apply Corollary \ref{cor_elemineq} under these conditions with
\begin{align*}
\epsilon_v(A)\lambda_1&=\log\Omega^{\epsilon_v(A)}_{\tilde D_{I_{1,v},v},v}(P) \geq \cdots \geq \epsilon_v(A)\lambda_{l_v}=\log\Omega^{\epsilon_v(A)}_{\tilde D_{I_{l_v,v},v},v}(P)\geq 0, \label{chain_ineq_ai}  
\end{align*}
\begin{align*}
b_1&=b_{1,v}-b_{0,v},\ldots, b_{l_v}=b_{l_v,v}-b_{l_v-1,v}, \nonumber 
\end{align*}
\begin{align*}
c_j&=c_{i_j,v}. \nonumber
\end{align*}
Using the telescoping property for the $b_i$,
\begin{equation}\label{cheb_application}
\log\prod_{v\in S} \prod_{s=1}^{l_v} \Omega^{\left(\max_{j} \frac{\sum_{j'=1} ^{j} c_{i_{j'},v}}{b_{j,v}}\right)(b_{s,v}-b_{s-1,v})\epsilon_v(A)}_{\tilde D_{I_{s,v},v},v}(P)
\geq\log\prod_{v\in S}\prod_{j=1}^{l_v} \Omega^{c_{i_j,v}\epsilon_v(A)}_{\tilde D_{I_{j,v},v},v}(P).
\end{equation}

It is known that the Seshadri constant of an intersection is bounded below by the minimum of the Seshadri constants of the components of the intersection, according to the explanation in terms of the reciprocal of the Seshadri constant, in \cite[Example 5.4.11]{PAGI}. Therefore, for any $j\in\{1,\ldots,l_v\}$, together with \eqref{sesh_bound}, we have
\begin{equation}
\epsilon_{\tilde D_{I_{j,v},v}}(A)\geq \min\{\epsilon_{\tilde D_{i_1,v}}(A),\ldots,\epsilon_{\tilde D_{i_j,v}}(A)\}\geq \epsilon_v(A)-\max_k \tau_{k,v} \cdot 2  \min_{l}  \epsilon_{\tilde  D_{l,v}}(A).\label{laz_ineq}
\end{equation}
Then \eqref{cheb_application} implies
\begin{align*}
& \log\prod_{v\in S}\prod_{s=1}^{l_v}\Omega^{\left(\max_{j} \frac{\sum_{j'=1} ^{j} c_{i_{j'},v}}{b_{j,v}}\right)(b_{s,v}-b_{s-1,v})\left(\epsilon_{\tilde{D}_{I_{s,v},v},v}(A)+\max_k \tau_{k,v} \cdot 2  \min_{l}  \epsilon_{\tilde  D_{l,v}}(A)\right)}_{\tilde D_{I_{s,v},v},v}(P) 
\end{align*}
\begin{align*}
\geq & \log\prod_{v\in S}\prod_{j=1}^{l_v} \Omega^{c_{i_j,v}\epsilon_v(A)}_{\tilde D_{I_{j,v},v},v}(P) \geq \log\prod_{v\in S}\prod_{j=1}^{l_v} \Omega^{c_{i_j,v}\epsilon_{\tilde D_{i_{j,v},v}}(A)}_{\tilde D_{I_{j,v},v},v}(P).
\end{align*}
If we construct a new list by repeating the divisors $\tilde{D}_{I_{s,v},v}$ exactly $b_{s,v} - b_{s-1,v}$ times (and omitting them if $b_{s,v} - b_{s-1,v} = 0$), then the resulting divisors will be in general position. By Theorem \ref{1.1}, there exists a proper Zariski closed subset $Z$ of $X$ such that
\begin{align*}
\log\prod_{v\in S} \prod_{s=1}^{l_v}\Omega^{(b_{s,v}-b_{s-1,v})\epsilon_{\tilde D_{I_{s,v},v}}(A) }_{\tilde D_{I_{s,v},v},v}(P)<(n+1+\epsilon')h_A(P)
\end{align*}
for all points $P\in X(k)\setminus Z$. We also find that by \eqref{modif_sesh} and additivity for local heights, for $i=1,\ldots,q$, we have
\begin{equation*}
\log\Omega^{\epsilon_{\tilde D_{i,v}}(A)}_{\tilde D_{i,v},v}(P)=\log\Omega^{\frac 1 {c_v(1+\tau_{i,v})} \cdot c_v(1+\tau_{i,v}) \epsilon_{D_{i,v}}(A)}_{D_{i,v},v}(P)= \log\Omega^{\epsilon_{D_{i,v}}(A)}_{D_{i,v},v}(P).
\end{equation*}
By selecting the quantities $\tau_{k,v}$ sufficiently small (relying solely on data related to the $D_{i,v}$, including finitely many domination inequalities of the form \eqref{lhf_ample_domination} for the $D_{i,v}$), we have, for a given $\epsilon''>0$
\begin{align*}
& \log\prod_{v\in S}  \prod_{s=1}^{l_v} \Omega^{(b_{s,v}-b_{s-1,v})(\max_k \tau_{k,v} \cdot 2  \min_{l}  \epsilon_{\tilde  D_{l,v}}(A))}_{\tilde D_{I_{s,v},v},v}(P)\\
\leq &\ \log\prod_{v\in S}  \prod_{s=1}^{l_v} \Omega^{(b_{s,v}-b_{s-1,v})(\max_k \tau_{k,v} \cdot 2 \epsilon_{\tilde  D_{i_{s,v},v}}(A))}_{\tilde D_{i_{s,v},v},v}(P)\\
\leq &\ \log\prod_{v\in S}  \prod_{s=1}^{l_v} \Omega^{(b_{s,v}-b_{s-1,v})(\max_k \tau_{k,v} \cdot 2 \epsilon_{D_{i_{s,v},v}}(A))}_{D_{i_{s,v},v},v}(P)\\
<&\ \epsilon''h_A(P).
\end{align*}
Taking $W=\Supp \tilde D_{I_{j,v},v}\neq \varnothing$, for $\alpha_v(W)$ as defined in Theorem \ref{1.7}
\begin{align*}
\frac{\sum_{j'=1} ^{j} c_{i_{j'},v}}{b_{j,v}}\leq \frac{\alpha_v(W)}{\codim W}.
\end{align*}
It should be noted that there are only finitely many possible combinations for indices $i_{1,v},\ldots, i_{l_v,v}$. Consequently, only finitely many applications of Theorem \ref{1.1} are required. Taking all of the above into consideration, there exists a proper Zariski-closed set $Z$ of $X$ such that
\begin{equation}\label{part1main}
\log\prod_{v\in S}\prod_{i=1}^q \Omega^{c_{i,v}\epsilon_{D_{i,v}}(A)}_{D_{i,v},v}(P)\leq \left((n+1)\max_{\substack{v\in S\\ \varnothing\subsetneq W\subsetneq X}}\frac{\alpha_v(W)}{\codim W}+\epsilon\right)h_A(P)
\end{equation}
for all points $P\in X(k)\setminus Z$.
		
The core of the next part of the proof is to apply iteratively the refined version of (\ref{part1main}). We first prove that (\ref{part1main}) holds for $D_{i,v}\sim d_{i,v}A$, with $c_{i,v}= 1$ for all $i$ and all $v\in S$.

Let $N$ be a positive integer such that $NA$ is very ample and $d_{i,v}$ divides $N$ for all $i$ and all $v\in S$. Define $M=l(NA)-1$. Let $\phi:X\to \mathbb{P}^M$ be the embedding of $X$ in $\mathbb{P}^M$ associated to $NA$. Since $\frac{N}{d_{i,v}}D_{i,v}\sim NA$, there exist hyperplanes $H_{i,v}$ in $\mathbb{P}^M$ defined over $k$ such that $\frac{N}{d_{i,v}}D_{i,v}=\phi^*H_{i,v}$. For $P\in X(k)\setminus \Supp D_{i,v}$, by the functoriality and additivity properties of Weil functions,
\begin{equation*}
\log\Omega_{H_{i,v},v}(\phi(P))=\log\Omega_{\frac{N}{d_{i,v}}D_{i,v},v}(P)=N\log\Omega^{\frac{1}{d_{i,v}}}_{D_{i,v},v}(P).
\end{equation*}

We also have $h(\phi(P))=Nh_A(P)+O(1)$ for $P\in X(k)$. Plugging these identities into (\ref{part1main}) for hyperplanes with $c_{i,v}=1$ for all $i$ and all $v\in S$, we have the desired result.
		
Ultimately, we demonstrate that the statement remains valid even if we replace linear equivalence with numerical equivalence. This conclusion is supported by the following result in algebraic geometry due to Matsusaka \cite[Theorem 1]{Mat} (see also \cite[Theorem 1.4.9]{PAGI} for details).

\begin{lemma}[Matsusaka]
\label{numample}
Let $A$ be an ample Cartier divisor on a projective variety $X$. For any Cartier divisor $D$ with $D\equiv NA$, there exists a positive integer $N_0$ such that for all $N\geq N_0$, $D$ is very ample.
\end{lemma}
The key point is to construct the following family of divisors $Q_{i,v}$ to replace $D_{i,v}$.

It is sufficient to prove the case where $d_{i,v}=1$ for all $i$ and all $v \in S$ by taking the least common multiple of the integers $d_{i,v}$, thus we have $D_{i,v}\equiv A$ for all $i$ and all $v \in S$. Let $N_0$ be the integer given by Lemma \ref{numample} for $A$. By standard properties of Weil functions and heights, there exists an integer $N>N_0$ sufficiently large such that
\begin{equation}
\label{Neq}
\log\prod_{v\in S}\prod_{i=1}^q\Omega^{\frac{N_0}{N}}_{D_{i,v},v}(P)< \epsilon h_A(P)+O(1) 
\end{equation}
for all $P\in X(k)\setminus \bigcup_{i,v}\Supp D_{i,v}$.

By the selection of $N_0$, it follows that $NA-(N-N_0)D_{i,v}$ is very ample for all $i$ and $v\in S$. According to this fact, there exist effective divisors $E_{i,v}$ such that $NA\sim (N-N_0)D_{i,v}+E_{i,v}$. Define $Q_{i,v}=(N-N_0)D_{i,v}+E_{i,v}$. For every subset \(I,J\subset\{1,\dots,q\}\) and every irreducible component
$$Z_I\subset \bigcap_{i\in I}\operatorname{Supp}D_{i,v}\cap X,$$
since the linear system is very ample, Bertini's theorem (see \cite[Chapter 3]{PAGI}) implies that the general choice of $E_{i,v}$ ensures
$$\dim\Bigl(Z_I\cap\bigcap_{j\in J}\operatorname{Supp}E_{j,v}\Bigr) \leq \dim Z_I-\sharp J.$$
Hence, for any irreducible closed subset $W\subset X$, defining
$I=\{i \mid W\subset \operatorname{Supp}D_{i,v}\}$ and $J=\{j \mid W\subset \operatorname{Supp}E_{j,v}\},$
we have
$$\dim W\le m-\sharp I-\sharp J.$$
Since $\alpha^{Q_{i,v}}_v(W)=\sharp\{i \mid W\subset \operatorname{Supp}Q_{i,v}\}\leq \sharp I+\sharp J$, it follows that $\alpha^{Q_{i,v}}_v(W)\le m-\dim W$, so the divisors $Q_{i,v}$
remain in $m$-subgeneral position. Moreover, the same estimate gives
$$\max_{\substack{v\in S\\ \varnothing\subsetneq W\subsetneq X}}\frac{\alpha^{Q_{i,v}}_v(W)}{\codim W}
= \max_{\substack{v\in S\\ \varnothing\subsetneq W\subsetneq X}}\frac{\sharp I}{\codim W}=\max_{\substack{v\in S\\ \varnothing\subsetneq W\subsetneq X}}\frac{\alpha_v(W)}{\codim W}.$$

Thus the subsequent application of (\ref{part1main}) to the $Q_{i,v}$ gives the same numerical constant as for the original $D_{i,v}$. For all $v\in S$, the divisors $(N-N_0)D_{1,v}+E_{1,v},\ldots, (N-N_0)D_{q,v}+E_{q,v}$ and $D_{1,v},\ldots,D_{q,v}$ have the same positional property. We now apply (\ref{part1main}) with $c_{i,v}=1$ to the linearly equivalent divisors $Q_{i,v}=(N-N_0)D_{i,v}+E_{i,v}$, $i=1,\ldots,q$, $v\in S$, along with $NA$.  
\begin{align*}
\log\prod_{v\in S}\prod_{i=1}^q \Omega_{(N-N_0)D_{i,v}+E_{i,v},v}(P)< \left((n+1)\max_{\substack{v\in S\\ \varnothing\subsetneq W\subsetneq X}}\frac{\alpha_v(W)}{\codim W}+\epsilon\right)h_{NA}(P)
\end{align*}
for all $P\in X(k)\setminus Z$ for some proper Zariski-closed subset $Z$ of $X$ containing the supports of all $E_{i,v}$, $v\in S$, $i=1,\ldots, q$, is established. Since $\log\Omega_{E_{i,v},v}$ is bounded from below outside the support of $E_{i,v}$, we obtain
\begin{equation*}
\log\prod_{v\in S}\prod_{i=1}^q\Omega^{\left(1-\frac{N_0}{N}\right)}_{D_{i,v},v}(P)< \left((n+1)\max_{\substack{v\in S\\ \varnothing\subsetneq W\subsetneq X}}\frac{\alpha_v(W)}{\codim W}+\epsilon\right)h_A(P)+O(1)
\end{equation*}
for all $P \in X(k)\setminus Z$. Inserting (\ref{Neq}) into the above formula, we deduce that
\begin{equation}
\label{eq1}
\log\prod_{v\in S}\prod_{i=1}^q\Omega_{D_{i,v},v}(P)<\left((n+1)\max_{\substack{v\in S\\ \varnothing\subsetneq W\subsetneq X}}\frac{\alpha_v(W)}{\codim W}+2\epsilon\right)h_A(P)+O(1) 
\end{equation}
for all $P\in X(k)\setminus Z$.

Given that there are only finitely many $k$-rational points with bounded height (with respect to $A$), we can ensure that equation (\ref{eq1}) holds without the $O(1)$ term by adding a finite number of points to $Z$.\par
		
$\underline{{Approach}-\mathbb{II}}$: Another version of the proof is obtained by using an immediate corollary of Theorem \ref{1.7}. When $D$ is an effective Cartier divisor, its Seshadri constant simplifies by the properties stated in Definition \ref{Sesh}. We recall that
\begin{align*}
\epsilon_D(A)=\sup\{\gamma\in {\mathbb{Q}}^{\geq 0}\mid A-\gamma D\text{ is $\mathbb{Q}$-nef}\},
\end{align*}
so that $\tilde{X}=X$. Therefore, we obtain
\begin{corollary}\label{generalized}
Let $X$ be a  projective variety of dimension $n$  defined over a number field $k$, and let $S$ be a finite set of places of $k$. For each $v\in S$, let $D_{1,v},\ldots, D_{q,v}$ be effective Cartier divisors of $X$, defined over $k$, and let $\gamma_{D_{i,v}} $ be rational numbers such that $A-\gamma_{D_{i,v}}D_{i,v}$ is a nef $\mathbb{Q}$-divisor for all $i$ and $v\in S$, where $A$ is an ample Cartier divisor. Let $\epsilon > 0$ and $c_{i,v} \in \mathbb{R}^{\geq 0}$. Then there exists a proper Zariski-closed subset $Z \subset X$ such that
\begin{align*}
\log\prod_{v\in S}\prod_{i=1}^q \Omega^{c_{i,v}\gamma_{D_{i,v}}}_{D_{i,v},v}(P)< \left((n+1)\max_{\substack{v\in S\\ \varnothing\subsetneq W\subsetneq X}}\frac{\alpha_v(W)}{\codim W}+\epsilon\right)h_A(P)
\end{align*}
holds for all points $P\in X(k)\setminus Z$.
\end{corollary} 

We directly set $c_{i,v}=1$ for all $i$ and $v\in S$, with the coefficient $\gamma_{D_{i,v}}=\frac{1}{d_{i,v}},$ hence we obtain (\ref{eq1}).
		
\subsection{Part II}\label{P2}

An inequality involving parameters and controlling the dimension of $Z$.\par
In this part, let $c_{i,v}=1$ for all $i$ and $v\in S$ and assume $\delta m \geq 2$.
In order to construct our results, we will verify the following key inequality that reflects the need for controlling exceptional-set contributions from all possible component dimensions $j \in \left\{l, \ldots , \delta n\right\}$.
\begin{align}\label{CBA}
C(\delta, l,m,n)\leq B(l,\delta m, \delta n)\leq \max_{l\leq j\leq \delta n}A(\delta m,j).
\end{align}
The purpose of proving this inequality is twofold.
First, we reformulate this theorem so that it takes the form of our results concerning Weil functions rather than in terms of the proximity functions defined in \cite[Section 2]{levin_duke}. Second, we reveal its underlying principles to refine our discussion in the next part and subsequent sections.

\begin{proof}
Consider the quantities and divisors $c , l , m , n , k , S , X , A , D_{1,v}, \ldots, D_{q,v}$, $q=m+1$ and $D$, as previously defined in Section \ref{R1.2}. According to Remark \ref{lem:seshadri_numerical}, we have $\epsilon_{D_{i,v}}=\frac{1}{d_{i,v}}$ for all $i$ and $v \in S$. 

Let $Z$ denote the exceptional set referred to in \eqref{exc1}. 
Assume that dim$Z=j \geq l$. For the right-hand side, we will only verify the case of rational points with $\delta=1$. The case $\delta>1$ will then follow from the symmetric power construction in the next step.
 
We restrict to an irreducible component $W$ of $Z$ (over $\kbar$) with $\dim W = j$. By extending the field $k$ if necessary, we may assume that $W$ is defined over $k$. The restricted divisors $D_{1,v}|_W, \ldots, D_{q,v}|_W$ are in $m$-subgeneral position and $A|_W$ remains ample on $W$. By the properties of Weil functions and the setup of $W$, the set of points $P\in W(k)\setminus \Supp D$ for which
\begin{align*}
\log\prod_{v\in S}\prod_{\substack{w\in M_{k(P)}\\w\mid v}} \prod_{i=1}^q\Omega^{\epsilon_{D_{{i,v}|W}}(A)}_{D_{{i,v}|W},w}(P)>(c-\epsilon) h_{A|_W}(P)
\end{align*}
is Zariski-dense in $W$. This shows that $A(m,j)\geq c$, which implies the right-hand side.

Now consider algebraic points of bounded degree $\delta$. For simplicity and without loss of generality, for each $v \in S$, we replace $D_{i,v}$ by $\left(\prod_{j\neq i}\frac{1}{\epsilon_{D_{j,v}}(A)}\right)D_{i,v}$ and the ample Cartier divisor $A$ by $\left(\prod_{j=1}^q\frac{1}{\epsilon_{D_{j,v}}(A)}\right)A$, so that $\epsilon_{D_{i,v}}(A)=1$ for all $i$ and all $v \in S$. We still write $D_{i,v}$ and $A$ for the normalized divisors, suppressing the dependence on $v$. 

For $[k(P):k]=\delta$ and $\epsilon>0$, define $Z$ to be the Zariski-closure of the set
\begin{align*}
\left\{P\in X(\kbar)\setminus \Supp D\mid  \log\prod_{v\in S}\prod_{\substack{w\in M_{k(P)}\\w\mid v}} \prod_{i=1}^q\Omega_{D_{i,v},w}(P)>(B(l,\delta m, \delta n)+\epsilon)h_A(P)\right\}.
\end{align*}
It suffices to show that $\dim Z<l$, which is equivalent to proving $C(\delta, l,m,n)\leq B(l,\delta m, \delta n).$

Now consider the following construction. 
			
Let $X^\delta$ denote the $\delta$-th Cartesian power of $X$. Let $X^{(\delta)}$ denote the $\delta$-th symmetric power of $X$, defined as the quotient $X^\delta / S_\delta$, where $S_\delta$ is the symmetric group on $\delta$ elements acting by permuting the coordinates of $X^\delta$. 
Let $$\phi: X^\delta \to X^{(\delta)}$$ be the natural map that takes a point in $X^\delta$ to its equivalence class.
			
For a point $P \in X(\kbar)$ satisfying $[k(P):k] = \delta$, let $P_1, P_2, \ldots, P_\delta$ denote the $\delta$ distinct conjugates of $P$ over $k$. Define the map $$\psi: X(\kbar) \to X^\delta$$ by $\psi(P) = (P_1, P_2, \ldots, P_\delta)$. One checks that the point $\phi(\psi(P)) \in X^{(\delta)}$ is $k$-rational.

Additionally, we define push-forward and pullback maps on effective Cartier divisors.			
Let $$\pi_i: X^\delta \to X$$ denote the $i$-th projection map, which projects a point in $X^\delta$ onto its $i$-th coordinate. Now we use these maps and $\phi$ to define the following divisors.
 
Define $$D^*_{i,v} = \phi_* \pi_1^* (D_{i,v}),$$ and the sum of all these divisors over all places in $S$, 
$$D^* = \sum_{v \in S} \sum_{i=1}^q D^*_{i,v}= \sum_{v \in S} \sum_{i=1}^q\phi_* \pi_1^* (D_{i,v}).$$ 
Also define $$A^* = \phi_* \pi_1^* (A).$$ 
Note that $$\phi^* \phi_* \pi_1^* D_{i,v} = \sum_{j=1}^\delta \pi_j^* D_{i,v}.$$
Since the $D_{i,v}$ are in $m$-subgeneral position, the $D^*_{i,v}$ are in $m\delta$-subgeneral position. It can be verified that $A^*$ is also ample. Consequently, the $D^*_{i,v}$ are effective and ample for all $i$ and $v$. If $$D_{i,v} \equiv D_{j,v} \equiv A$$ for all $i$ and $j$ , then $$D^*_{i,v} \equiv D^*_{j,v} \equiv A^*$$ as well.
			
Thus by the definition of $B(l,\delta m, \delta n)$, we have the inequality
\begin{equation}
\label{beq}
\log\prod_{v\in S}\prod_{\substack{w\in M_{k(P)}\\w\mid v}} \prod_{i=1}^q\Omega_{D^*_{i,v},w}(P)<(B(l,\delta m, \delta n)+\epsilon)h_{A^*}(P)
\end{equation}
\begin{align*}
\quad \forall P\in X^{(\delta)}(k)\setminus (Z'\cup \Supp D^*),
\end{align*}
for some Zariski-closed subset $Z'$ of $X^{(\delta)}$ satisfying $$\dim Z'\leq l-1.$$
For $P\in X(\kbar)\setminus \Supp D$, $[k(P):k]=\delta$, we also have the following chain of equalities, up to $O(1)$,
\begin{align*}
&\log\prod_{v\in S}\prod_{\substack{w\in M_{k(P)}\\w\mid v}} \prod_{i=1}^q\Omega_{D^*_{i,v},w}(\phi(\psi(P)))=
\end{align*}
\begin{align*}
&\log\prod_{v\in S}\prod_{\substack{w\in M_{k(P)}\\w\mid v}} \prod_{i=1}^q\Omega_{\phi_*\pi_1^*D_{i,v},w}(\phi(\psi(P)))=
\end{align*}
\begin{align*}
\log\prod_{v\in S}\prod_{\substack{w\in M_{k(P)}\\w\mid v}} \prod_{i=1}^q\Omega_{\phi^*\phi_*\pi_1^*D_{i,v},w}(\psi(P))
=\log\prod_{v\in S}\prod_{\substack{w\in M_{k(P)}\\w\mid v}} \prod_{i=1}^q\Omega_{\sum_{j=1}^\delta \pi_j^*D_{i,v},w}(\psi(P))
\end{align*}
			
\begin{align*}
=&\sum_{j=1}^{\delta}\log\prod_{v\in S}\prod_{\substack{w\in M_{k(P)}\\w\mid v}} \prod_{i=1}^q\Omega_{\pi_j^*D_{i,v},w}(\psi(P))=
\end{align*}
\begin{align*}
&\sum_{j=1}^{\delta}\log\prod_{v\in S}\prod_{\substack{w\in M_{k(P)}\\w\mid v}} \prod_{i=1}^q\Omega_{D_{i,v},w}(\pi_j\psi(P))=\sum_{j=1}^{\delta}\log\prod_{v\in S}\prod_{\substack{w\in M_{k(P)}\\w\mid v}} \prod_{i=1}^q\Omega_{D_{i,v},w}(P_j)
\end{align*}
\begin{align*}
=\delta\log\prod_{v\in S}\prod_{\substack{w\in M_{k(P)}\\w\mid v}} \prod_{i=1}^q\Omega_{D_{i,v},w}(P).
\end{align*}

Similarly, we also have $h_A^*(\phi(\psi(P)))=\delta h_A(P)+O(1)$. Consider $Z^*\subset Z(\kbar)$ as a Zariski-dense set of points in $Z$ satisfying $[k(P):k]=\delta$. It follows that for an appropriate choice of the height functions,
\begin{equation}\label{refined divisor}
\log\prod_{v\in S}\prod_{\substack{w\in M_{k(P)}\\w\mid v}} \prod_{i=1}^q\Omega_{D^*_{i,v},w}(\phi(\psi(P)))>(B(l,\delta m, \delta n)+\epsilon)h_{A^*}(\phi(\psi(P)))
\end{equation}
for all $P\in Z^*$. From \eqref{beq}, we have $\dim \phi(\psi(Z^*))\leq l-1$.  Since $\phi$ is a finite map and $\pi_1(\psi(Z^*))=Z^*$, $\dim \phi(\psi(Z^*))\geq \dim Z^*=\dim Z$.  Thus $\dim Z\leq l-1$, as desired. This shows that $C(\delta, l,m,n)\leq B(l,\delta m, \delta n)$. Based on the above construction and discussion, we can also obtain $B(l, \delta m, \delta n)\leq \max_{l\leq j\leq  \delta n}A(\delta m,j)$.
\end{proof}
		
\subsection{Part III}\label{P3}

Part II for algebraic points of bounded degree.

In this part, let $m$ be a positive integer with $m \geq n$ and set $q:=m+1$.
		
Using formula (\ref{eq1}) in Part I and combining it with inequality (\ref{CBA}) in Part II, we can immediately obtain Theorem \ref{main} with the factor of $a_\delta+\epsilon$. Now we introduce another factor. 

According to the fact that
\begin{align}\label{two distri}
	\delta_{\mathcal{D}_v}=\max_{I\subset \{1,\ldots, q\}}\max\left\{1, \frac{\sharp I}{\codim \bigcap_{i\in I}\Supp D_{i,v}}\right\}=\max\left\{1, \max_{\varnothing\subsetneq W\subsetneq X}\frac{\alpha_v(W)}{\codim W}\right\},
\end{align}
we can replace $a_\delta$ on the right side by a slightly weakened factor $a_{\Delta,\delta}$,
where $a_{\Delta,\delta}$ here is defined by
\begin{align*}
\max_{1\leq j\leq \delta n}\iota(\delta m,j),
\end{align*}
which is related to the distributive constant. Define $\iota(m,n)=(m-n+1)(n+1)$ for $m$-subgeneral position or $(\frac{m-n}{\kappa}+1)(n+1)$ for $m$-subgeneral position with index $\kappa$.
 
To explain this, we show here that the above two definitions in (\ref{two distri}) are equivalent.
\begin{proof}		
Assume that $\tilde{W}$ is a subvariety of $X$ at which the ratio $$\frac{\alpha_v(W)}{\codim W}$$ attains its maximum. 

By reordering if necessary, we can assume that $\tilde{W} \subset \Supp D_{i,v}$ for $i = 1, \dots, \alpha_v(\tilde{W})$. Define $\bar{W} = \bigcap_{i=1}^{\alpha_v(\tilde{W})} \Supp D_{i,v}$. It is evident that $\tilde{W} \subset \bar{W}$, which implies that $\codim \tilde{W} \geq \codim \bar{W}$. Additionally, since $\tilde{W} \nsubseteq \Supp D_{i,v}$ for all $i > \alpha_v(\tilde{W})$, it follows that $\bar{W} \nsubseteq \Supp D_{i,v}$ for all $i > \alpha_v(\tilde{W})$. Therefore, we conclude that $\alpha_v(\tilde{W})=\alpha_v(\bar{W})$. Thus, we deduce that
$$\frac{\alpha_v(\bar{W})}{\codim \bar{W}} \geq \frac{\alpha_v(\tilde{W})}{\codim \tilde{W}}.$$
By the assumption on $\tilde{W}$,
$$\frac{\alpha_v(\bar{W})}{\codim \bar{W}}=\frac{\alpha_v(\tilde{W})}{\codim \tilde{W}}.$$
Therefore, it suffices to consider those $W$ that are the intersections of some elements in $\{D_{i,v}\}$.
\end{proof}

At this point, in conjunction with the inequality and assumptions in $\S$\ref{P2}, we obtain Proposition \ref{pp}. The inequality 
\begin{align*}
\log\prod_{v\in S}\prod_{\substack{w\in M_{k(P)}\\w\mid v}}\prod_{i=1}^q \Omega^{\epsilon_{D_{i,v}}(A)}_{D_{i,v},w}(P)< \left(a_{\Delta,\delta}+\epsilon\right)h_A(P)+O(1)
\end{align*}
holds for all points $P\in X(\kbar)\setminus \bigcup_{i,v}\Supp D_{i,v}$, where the $O(1)$-term depends on $\epsilon$.
		
Compared with Theorem \ref{main}, because of the particular properties of the other cases, we require a unified and  explicit representation that is broader in scope, which motivates the form of our results involving $a_{\Delta,\delta}$. Moreover, geometrically, this representation incorporates the contributions arising from different positional conditions.
		
If $D_{i,v}$ are taken to be effective Cartier divisors in $m$-subgeneral position, then 
\begin{align*}
\max_{v\in S}\delta_{\mathcal{D}_v}(n+1)=\Delta(n+1)\leq (m-n+1)(n+1)=\iota(m,n).
\end{align*}
The inequality (\ref{CBA}) shows that
\begin{align*}
C(\delta, 1,m,n)\leq B(1,\delta m, \delta n)\leq \max_{1\leq j\leq \delta n}A(\delta m,j)
\end{align*}
\begin{align*}
\leq \max_{1\leq j\leq \delta n}\left(\delta m-j+1\right)(j+1)=\max_{1\leq j\leq \delta n}\iota(\delta m,j).
\end{align*}
Viewing $$\left(\delta m-j+1\right)(j+1), \quad 1\leq j\leq \delta n$$ as a quasi-quadratic function of the integer variable $j$, we divide the discussion into two parts (the definition of $[*]$ can be found in Section \ref{a2}).
 
(i) For $m>2n$,
$$a_{\Delta,\delta}={\delta^2}n(m-n)+{\delta}m+1.$$

(ii) For $n \leq m \leq 2n$, setting $i=\left[\frac{\delta m}{2}\right]$, we obtain the maximum $$a_{\Delta,\delta}=-i^2+\delta m \cdot i+{\delta}m+1. $$ 

Now we show that $a_{\Delta,\delta}$ for $n \leq m \leq 2n$ is also better than the result in Theorem \ref{1.9}. 
\begin{proof}
We consider the parity of $\delta m$:
Let
$$
x=\delta m,\qquad y=\delta n,\qquad i=\left[\frac{x}{2}\right].
$$
Since $n \leq m \leq 2n$, we have $y\leq x\leq 2y$. Assume $x>2$.

Define $$ A=-i^2+xi+x+1,\quad B=\frac{x(x-1)(y+1)}{x+y-2}. $$

To compare $A$ and $B$, consider
$$ N=x(x-1)(y+1)-A(x+y-2).$$
Since the denominator $x+y-2>0$, the sign of $B-A$ is the same as the sign of $N$.

\underline{Case $\mathbb{I}$: $x=2k$ is even.}

Since $x>2$, we have $k\geq 2$. Then
$$
i=k,\qquad A=(k+1)^2.
$$
A calculation shows
$$
N=(3k^2-4k-1)y-2k^3+2k^2+2.
$$
The coefficient of $y$ is $3k^2-4k-1>0$, so $N$ is increasing in $y$. Since
$$
y\ge \frac x2=k,
$$
the minimum occurs at $y=k$. Substituting $y=k$,
$$ N_{\min}=(k-1)(k-2)(k+1).$$
\begin{itemize}
	\item If $k=2$, we have $x=4$, then $N_{\min}=0$, and equality occurs only when $y=k=2$, i.e.
	$$
	\delta m=4,\quad \delta n=2.
	$$
	\item If $k>2$, we have $x\geq 6$, then $N_{\min}>0$, so $B>A$.
\end{itemize}

\underline{Case $\mathbb{II}$: $x=2k+1$ is odd.}

Since $x>2$, we have $k\geq 1$. Then $i=k$ or $i=k+1$, and
$$A=(k+1)(k+2).$$
A calculation gives
$$ N=(3k^2-k-2)y-2k^3-k^2+k+2.$$
The coefficient is
$$ 3k^2-k-2=(3k+2)(k-1).$$
\begin{itemize}
	\item If $k=1$, we have $x=3$, then both the coefficient and the constant term are $0$, so $N=0$. Hence $A=B$.
	\item If $k\geq 2$, we have $x\geq 5$, then the coefficient is positive. Since
	$$ y\ge k+1,$$
	the minimum occurs at $y=k+1$. Substituting $y=k+1$,
	$$ N_{\min}=k(k-1)(k+2)>0.$$
	Hence $B>A$.
\end{itemize}

Therefore, when $\delta m>2,$ we always have
$$B\ge A,$$
that is,
$$ \frac{\delta m(\delta m-1)(\delta n+1)}{\delta m+\delta n-2}\ge -i^2+\delta m\,i+\delta m+1.$$
Equality holds if and only if
$$ \delta m=3 $$
or
$$\delta m=4,\quad \delta n=2.$$
\end{proof}
And if $\mathcal{D}_v=\left\{D_{i,v}\right\}$ is a set of effective Cartier divisors in $m$-subgeneral position with index $\kappa$, which is already in $m$-subgeneral position according to Definition \ref{1.5}, then
\begin{align*}
\Delta\leq\frac{m-n}{\kappa}+1, \quad \iota(m,n)= (\frac{m-n}{\kappa}+1)(n+1),
\end{align*}
which implies that
\begin{align*}
C(\delta, 1,m,n)\leq B(1,\delta m, \delta n)\leq \max_{1\leq j\leq \delta n}A(\delta m,j)\leq \max_{1\leq j\leq \delta n}\left(\frac{\delta m-j}{\kappa}+1\right)(j+1)=\max_{1\leq j\leq \delta n}\iota(\delta m,j).
\end{align*}
By the same operation as above, we assert that the following factor estimates hold.

(i) For $m \geq \lceil 2n+\frac{1-\kappa}{\delta} \rceil$,
$$a_{\Delta,\delta}=\frac{n}{\kappa}\left({\delta^2}m-{\delta^2}n+\delta(\kappa-1)\right)+\frac{{\delta}m}{\kappa}+1.$$
 
(ii) For $n \leq m \leq  \lfloor 2n+\frac{1-\kappa}{\delta} \rfloor $, we have $i=\left[\frac{\delta m+(\kappa-1)}{2}\right]$ in this situation. Thus
$$a_{\Delta,\delta}=-\frac{i^2}{\kappa}+\left(\frac{\delta m+(\kappa -1)}{\kappa}\right)i+\frac{{\delta}m}{\kappa}+1.$$

This completes the proof.

\section{Supplementary propositions}
\label{sp}
As a powerful tool for measuring positivity, Seshadri constants relax stringent conditions in a certain sense. This yields a version that allows us to weaken the nefness assumption in our theorems. With some modifications, we can obtain the following proposition. It is worth noting that we retain the original subgeneral position assumption while introducing weights. This broadens the scope of applicability of the following generalization, although it also further weakens our results.

\begin{proposition}\label{generalized-ap}
Let $X$ be a  projective variety of dimension $n$  defined over a number field $k$, and let $S$ be a finite set of places of $k$. Let $m \geq n$ be a positive integer and set $q=m+1$. For each $v\in S$, let $D_{1,v},\ldots, D_{q,v}$ be effective Cartier divisors on $X$, defined over $k$ and in $m$-subgeneral position. Let  $\gamma_{D_{i,v}} $ be rational numbers such that $A-\gamma_{D_{i,v}}D_{i,v}$ is a nef $\mathbb{Q}$-divisor for all $i$ and for all $v\in S$ , where $A$ is an ample Cartier divisor on $X$. Let $c_{1,v},\ldots, c_{q,v}$ be nonnegative real numbers and let $\delta$ be a positive integer. For a closed subset $W\subset X$ and $v \in S$, let $a_\delta$ be as defined in Section \ref{a1} with
\begin{align*}
	\alpha_v(W)=\sum_{\substack{i\\ W\subset \Supp D_{i,v}}}c_{i,v}.
\end{align*}
Then for every $\epsilon > 0$, there exists a proper Zariski-closed subset $Z \subset X$ such that
\begin{align*}
\log\prod_{v\in S}\prod_{\substack{w\in M_{k(P)}\\w\mid v}}\prod_{i=1}^q \Omega^{c_{i,v}\gamma_{D_{i,v}}}_{D_{i,v},w}(P)< \left(a_{\delta}+\epsilon\right)h_A(P)+O(1)
\end{align*}
for all points $P\in X(\kbar)\setminus \bigcup_{i,v}\Supp D_{i,v}$ satisfying $[k(P):k] \leq \delta$, where the $O(1)$-term depends on $\epsilon$.
\end{proposition}

\begin{proof}
Assume $q = m + 1$.

Let $\gamma_{D_{i,v}} $ be rational numbers such that $$A-\gamma_{D_{i,v}}D_{i,v}$$ is a nef $\mathbb{Q}$-divisor for all $i$ and for all $v\in S$. 
	
Consider the quantities and divisors in $m$-subgeneral position $c , l , m , n , k , S , X , A , D_{1,v}, \ldots, D_{q,v}$, and $D$, as previously defined. For $\delta\in \mathbb{Z}^{+}$, define the geometric exceptional set $$Z=Z\left(c,\delta, m, n,k, S, X, A, \{D_{i,v}\}\right)$$ to be the smallest Zariski-closed subset of $X$ such that the set
\begin{equation*}
\left\{P\in X(\kbar)\setminus \Supp D\mid \log\prod_{v\in S}\prod_{\substack{w\in M_{k(P)}\\w\mid v}} \prod_{i=1}^q\Omega_{D_{i,v},w}^{c_{i,v}\gamma_{D_{i,v}}}(P)>c h_A(P)\right\}\setminus Z
\end{equation*}
satisfying $[k(P):k]\leq \delta$ is finite for all choices of the Weil and height functions. Analogous definitions to those in Section \ref{R1.2} hold for $C(\delta, l,m,n)$, $B(l,\delta m, \delta n)$ and $A(\delta m,j)$. 
		
If, for rational points with $\delta=1$, we use the same argument and hypothesis as in Part II, we have that
\begin{align*}
\log\prod_{v\in S}\prod_{\substack{w\in M_{k(P)}\\w\mid v}} \prod_{i=1}^q\Omega^{c_{i,v}\gamma_{D_{{i,v}|W}}}_{D_{{i,v}|W},w}(P)>(c-\epsilon) h_{A|_W}(P)
\end{align*}
is Zariski-dense in $W$. This implies the right-hand side of the inequality (\ref{CBA}) for rational points with $\delta=1$ under this nef condition (this also holds true for algebraic points with $\delta > 1$).
		
Let $\epsilon>0$ and define $Z$ to be the Zariski-closure of the set
\begin{multline*}
\{P\in X(\kbar)\setminus \Supp D\mid [k(P):k]=\delta, \\  
\log\prod_{v\in S}\prod_{\substack{w\in M_{k(P)}\\w\mid v}} \prod_{i=1}^q\Omega^{c_{i,v}\gamma_{D_{i,v}}}_{D_{i,v},w}(P)>(B(l,\delta m, \delta n)+\epsilon)h_A(P)\}
\end{multline*}

To simplify the proof, and without loss of generality, we replace $D_{i,v}$ by $\left(\prod_{j\neq i}\frac{1}{\gamma_{D_{j,v}}}\right)D_{i,v}$ and $A$ by $\left(\prod_{j=1}^q\frac{1}{\gamma_{D_{j,v}}}\right)A$ for each $v \in S$. That is, $\gamma_{D_{i,v}}=1$ for all $i$ and $v \in S$.

With the same mapping construction as in $\S$\ref{P2}, we set $$D^*_{i,v} = \phi_* \pi_1^* (D_{i,v}), \quad A^* = \phi_* \pi_1^* (A).$$ 
The key issue is to verify that such a mapping combination also makes each divisor $$A^* -\gamma_{D_{i,v}}D_{i,v}^*$$ $\mathbb{Q}$-nef. By construction,
$$\phi^*(A^*-\gamma D_{i,v}^*)=\sum_{j=1}^{\delta}\pi_j^*(A-\gamma D_{i,v}).$$
Since $A-\gamma D_{i,v}$ is $\mathbb{Q}$-nef, we have that each
$\pi_j^*(A-\gamma D_{i,v})$ is nef, and their finite sum is also nef. Thus
$$\phi^*(A^*-\gamma D_{i,v}^*)\ \text{is $\mathbb{Q}$-nef.}$$

It is known that the nefness descends along finite surjective morphisms: if $L$ is a
$\mathbb{Q}$-Cartier divisor on $X^{(\delta)}$ and $\phi^*L$ is nef,
then $L$ is nef. Indeed, for every curve $C\subset X^{(\delta)}$, since $\phi^*C$ is an effective $1$-cycle and $\deg(\phi)>0$, we have
$$\phi^*L\cdot \phi^*C
=\deg(\phi)\,L\cdot C \geq 0.$$
Therefore,
$$A^*-\gamma D_{i,v}^*\ \text{is }\mathbb{Q}\text{-nef}.$$
 
The remaining steps follow from (\ref{beq}) and (\ref{refined divisor}). Combining this with Corollary \ref{generalized} and the definition of $a_\delta$, we obtain the result.
\end{proof}
	
Invoking \cite[Theorem 2.5.8]{Vojta_LNM}, where $I$ runs through all subsets of $\left\{1,\dots, q\right\}$ ($q \geq m+1$) such that the divisors $D_{i}$, $i \in I$, are in $m$-subgeneral position, we can extend Proposition \ref{generalized-ap} to
\begin{proposition}\label{p4.2}
Let $X$ be a  projective variety of dimension $n$ defined over a number field $k$, and let $S$ be a finite set of places of $k$. Let $m \geq n$ be a positive integer, and let $q \geq m+1$ be an integer. Let $D_{1},\ldots, D_{q}$ be effective Cartier divisors on $X$, defined over $k$. Let $\gamma_{D_{i}}$ be a rational number such that $A-\gamma_{D_{i}}D_{i}$ is a nef $\mathbb{Q}$-divisor for all $i$, where $A$ is an ample Cartier divisor on $X$. Let $c_{i} \in \mathbb{R}^{\ge 0}$ for all $i$, and let $\delta$ be a positive integer. For a closed subset $W \subset X$, define $$\alpha(W)=\sum_{\substack{i\\ W\subset \Supp D_{i}}}c_{i}, \quad \ell (m,n)=(n+1)\max_{\substack{\varnothing\subsetneq W\subsetneq X}}\frac{\alpha(W)}{\codim W}.$$ Then, we have the inequality
\begin{align*}
\log\prod_{v\in S}\prod_{\substack{w\in M_{k(P)}\\w\mid v}}\max_I\prod_{i \in I} \Omega^{c_i\gamma_{D_{i}}}_{D_{i,w}}(P)< \left(a_{\delta}+\epsilon\right)h_A(P)+O(1)
\end{align*} 
for all points $P\in X(\kbar)\setminus \bigcup_{i}\Supp D_{i}$ satisfying $[k(P):k] \leq \delta$, where the $O(1)$-term depends on $\epsilon$.
\end{proposition}

\begin{proof} 
		
Since the divisors $D_i$ are in $m$-subgeneral position, any point $P\in X(k)$ can be $v$-adically close to at most $m$ of the divisors $D_i$.  This implies that there exists a constant $c$ such that for any $P\in X(k)\setminus \bigcup_{i=1}^q\Supp D_i$ and any $v\in S$, there are indices $i_1,\ldots, i_{m'}\in I \subset \{1,\ldots, q\}$, $m' \leq m$ such that
\begin{equation*}
\log\prod_{i \in I} \Omega_{D_i,v}^{c_i\gamma_{D_{i}}}(P)<\log\prod_{j=1}^{m'}\Omega_{D_{i_j},v}^{c_{i_j}\gamma_{D_{i_j}}}(P)+c.
\end{equation*}
After reindexing the divisors $D_{i}$, it suffices to prove a slightly weaker inequality involving at most $m$ divisors instead of $m+1$ as in Proposition \ref{generalized-ap}. Note that there are only finitely many possibilities for the subset $I$ and for the divisors $\left\{D_{i_{m'}}\right\}$. Then, applying Corollary \ref{generalized} finitely many times with the factor  $$(n+1)\max_{\substack{\varnothing\subsetneq W\subsetneq X}}\frac{\alpha(W)}{\codim W},$$ we conclude that there exists a proper Zariski-closed subset $Z \subset X$ and $\epsilon > 0$ such that
\begin{align*}
\log\prod_{v\in S}\max_I\prod_{i \in I} \Omega^{c_i\gamma_{D_{i}}}_{D_{i},v}(P)< \left((n+1)\max_{\substack{\varnothing\subsetneq W\subsetneq X}}\frac{\alpha(W)}{\codim W}+\epsilon\right)h_A(P)+O(1).
\end{align*}
Combining this with the reformulation of the result in $\S$\ref{P2} implies our proposition immediately, which is a weaker but more general statement.
\end{proof}
\begin{remark}\label{4.3}
If $\mathcal{D}=\{D_{1},\ldots, D_{q}\}$ is assumed to be in $m$-subgeneral position, then $I$ includes the subset $\left\{1,\dots, q\right\}$, thus we obtain that for $q \geq m+1$, the inequality 
\begin{align*}
	\log\prod_{v\in S}\prod_{\substack{w\in M_{k(P)}\\w\mid v}}\prod_{i=1}^q \Omega^{c_{i}\gamma_{D_{i}}}_{D_{i},w}(P)< \left(a_{\delta}+\epsilon\right)h_A(P)+O(1)
\end{align*}
holds for all points $P\in X(\kbar)\setminus \bigcup_{i}\Supp D_{i}$ satisfying $[k(P):k] \leq \delta$.

Also, for our purposes, if we set $c_{i}=1$ for all $i$ and take $D_i$ ($1 \leq i \leq q$) to be hyperplanes, then we have statement (\ref{4}), with $a_\delta$ on the right-hand side.
\end{remark}
Taking $c_{i}=1$ for all $i$, we also have the following corollary, which incorporates the factor $a_{\Delta, \delta}$. We can thus omit the $m$-subgeneral position assumption.
\begin{corollary}\label{c4.4}
Let $X$ be a projective variety of dimension $n$ defined over a number field $k$, and let $S$ be a finite set of places of $k$. Fix a positive integer $m \geq n$ and an integer $q \geq m+1$. Let $D_{1},\ldots, D_{q}$ be effective Cartier divisors on $X$, defined over $k$. For each $i$, let $\gamma_{D_{i}}$ ($i=1,\dots,q$) be a rational number such that $A-\gamma_{D_{i}}D_{i}$ is a nef $\mathbb{Q}$-divisor, where $A$ is an ample Cartier divisor on $X$. Let $\delta$ be a positive integer, and let $a_{\Delta,\delta}$ be as defined in Section \ref{a2}.

Then for every $\epsilon > 0$, there exists a proper Zariski-closed subset $Z \subset X$ such that
	\begin{align*}
		\log\prod_{v\in S}\prod_{\substack{w\in M_{k(P)}\\w\mid v}}\prod_{i=1}^{q} \Omega^{\gamma_{D_{i}}}_{D_{i,w}}(P)< \left(a_{\Delta,\delta}+\epsilon\right)h_A(P)+O(1)
	\end{align*} 
	for all points $P\in X(\kbar)\setminus \bigcup_{i}\Supp D_{i}$ satisfying $[k(P):k] \leq \delta$, where the $O(1)$-term depends on $\epsilon$.
\end{corollary}

\section{An optimal factor and a conjecture}\label{OC}
	
It is notable that once we set $c_{i}=1$ for all $i$, we can construct the following exceptional case based on \cite[Example 1.3]{HL25}, which demonstrates the sharpness of $a_\delta$ and provides a crucial illustration for a proposition on hyperplanes.
	
Let $m$ hyperplanes intersect at a point $P_0$ but otherwise intersect generally. Additionally, for a positive integer $k$, choose a family of distinct hyperplanes $\left\{D_i\right\}$ ($i=1, \ldots, k$), each of which intersects the above $m$-hyperplanes generally and intersects generally in pairs. For $i=1, \ldots, k$ and $s=1,\ldots,\codim D_i$, set $D_{m+i_s}=D_i$. Set $q:=m+\sum_{i=1}^{k}\codim D_i$, and let $\mathcal{D}:=\left\{D_1, \dots, D_q\right\}$ be the resulting family of these $q$ hyperplanes. The family $\mathcal{D}$ is thus in $m$-subgeneral position. 
	
According to the initial conditions, there exists a maximum as follows.
$$\max_{\varnothing\subsetneq W\subsetneq X}\frac{\alpha(W)}{\codim W}=\frac{m}{n},$$
which, by the results in the proof section, gives the following formula
$$a_\delta=\max_{1\leq j\leq \delta n}\ell(\delta m,j)=\max_{1\leq j\leq \delta n}\left(\frac{\delta m}{j}\right)(j+1)=2\delta m.$$
Thus Proposition \ref{p4.2} for hyperplanes yields the following inequality
\begin{align*}
\log\prod_{v\in S}\prod_{\substack{w\in M_{k(P)}\\w\mid v}}\prod_{i=1}^q\left(\frac{\|f_{i}\|_{w}\max_j \|x_{j}\|_{w}}{\|f_{i}(P)\|_{w}}\right)< \left(2\delta m+\epsilon\right)h(P)+O(1)
\end{align*}
for all points $P\in X(\kbar)\setminus \Supp\sum_{i=1}^{q} D_i$, where the $O(1)$-term depends on $\epsilon$. 
	
\begin{remark}
Consider a line $L$ passing through $P_0$ and not contained in any of the hyperplanes $D_1, \ldots, D_m$. Such a line $L$ intersects $D_1, \ldots, D_q$ in at most two distinct points and contains an infinite Zariski-dense set of $(\sum_{i=1}^{q}D_i,S)$-integral points. Consequently, the factor $2\delta m+\epsilon$ cannot be replaced by any smaller value.
\end{remark}
If $m=n \geq 2$ (so that the position is general), then we certainly have 
$$C(\delta,1,n,n) \leq B(1,\delta n, \delta n) \leq \max_{1\leq j\leq \delta n}A(\delta n,j) \leq \max_{1\leq j\leq \delta n}\left(\frac{\delta n}{j}\right)(j+1)=2\delta n.$$

Now we state a proposition based on the above inequality.
\begin{proposition}\label{optimal factor} 
With notation as above, we assert that for $X=\mathbb{P}^n=\mathbb{P}^N$ and hyperplanes $\{D_i\}$,  $$C(\delta,1,n,n)=B(1,\delta n, \delta n)=2\delta n.$$ 
\end{proposition}
	
To verify this statement, we also need the following example due to Levin,
\begin{example}[{\cite[Example 7.4]{levin_duke}}]
\label{EX}
Let $\delta$ and $n\geq 2$ be positive integers. Consider $D_1,\ldots, D_{2\delta n}\subset \mathbb{P}^n$ as hyperplanes in general position, defined by homogeneous polynomials $f_{i}\in k[x_0,\ldots,x_n], 1 \leq i \leq 2\delta n$ over a number field $k$ such that the intersection $\bigcap_{i=(s-1)n+1}^{sn}D_i=\{P_s\}$ consists of a single point for $s=1,\ldots,2\delta$, where $P_1,\ldots, P_{2\delta}$ all lie on a line $L$ over $k$. It can be verified that there exists a finite set of places $S$ of $k$ and an infinite set of points $R\subset L(\kbar)\subset \mathbb{P}^n(\kbar)$ such that 
\begin{equation*}
\log\prod_{v\in S}\prod_{\substack{w\in M_{k(P)}\\w\mid v}}\prod_{i=1}^{2\delta n}\left(\frac{\|f_{i}\|_{w}\max_j \|x_{j}\|_{w}}{\|f_{i}(P)\|_{w}}\right)=2\delta nh(P)+O(1) 
\end{equation*}
and $[k(P):k]\leq \delta$ for all $P$ in $R$.
\end{example}
	
This example implies that $C(\delta,1,n,n)\geq2\delta n.$ 
Hence, we have Proposition \ref{optimal factor}. 

When the hyperplanes $\left\{D_{i}\right\}$ are in general position, this new bound $2\delta n$ for $C(\delta,1,n,n)$, which we achieved for the first time, not only substantially improves the factor compared with the previous one
$$(\delta n)^2\left(\frac{\delta n-1}{2\delta n-3}\right),$$ but also completely settles the following open problem in \cite{RVC} in the hyperplane case.

\begin{conjecture}[Conjecture \ref{c}]\label{cA}
Let $k$ be a number field, let $M_k$ be the set consisting of all the
non-equivalent places of $k$. Let $S \subset M_k$ be a finite set containing all the archimedean places. Let $\{L_{1},\ldots, L_{q}\}$ be linear forms with coefficients in $k$, in general position. Let $\delta$ be a positive integer such that $[k(\mathbf{x}):k]\leq \delta$. If $\epsilon > 0$, and if $C$ is a positive constant, then the set of points $\mathbf{x} \in [x_0:x_1:\cdots:x_n]$ with $x_j \in \kbar, 0 \leq j \leq n$, such that 
\begin{align}\label{cj}
\log\prod_{v\in S}\prod_{\substack{w\in M_{k(\mathbf{x})}\\w\mid v}}\prod_{i=1}^q \left(\frac{\|L_{i}\|_{w}\max_j \|x_{j}\|_{w}}{\|L_{i}(\mathbf{x})\|_{w}}\right)\geq\left(2\delta n+\epsilon\right)h(\mathbf{x})+C
\end{align}
is contained in a finite union of linear subspaces of dimension $0$.

We now prove that Proposition \ref{optimal factor} implies Conjecture \ref{cA}.
\begin{proof}
Let $\epsilon>0$ and $C>0$ be given, and set $c:=2\delta n+\epsilon/2$. Since $c>C(\delta,1,n,n)=2\delta n$, the definition of $C(\delta,1,n,n)$ yields a Zariski-closed set $Z_c\subset\mathbb P^n$ with $\dim Z_c<1$ (hence $Z_c$ is a finite set) such that the set of points $\mathbf{x}$
with $[k(\mathbf{x}):k]\leq\delta$ satisfying
\[
\log\prod_{v\in S}\prod_{\substack{w\in M_{k(\mathbf{x})}\\w\mid v}}\prod_{i=1}^q \left(\frac{\|L_{i}\|_{w}\max_j \|x_{j}\|_{w}}{\|L_{i}(\mathbf{x})\|_{w}}\right)>c\,h(\mathbf{x}).\]
is contained in $Z_c\cup F$ for some finite set $F$. Now for any $\mathbf{x}$ satisfying the inequality (\ref{cj}), we have
\[
\log\prod_{v\in S}\prod_{\substack{w\in M_{k(\mathbf{x})}\\w\mid v}}\prod_{i=1}^q \left(\frac{\|L_{i}\|_{w}\max_j \|x_{j}\|_{w}}{\|L_{i}(\mathbf{x})\|_{w}}\right)\geq(2\delta n+\epsilon)h(\mathbf{x})+C>(2\delta n+\epsilon/2)h(\mathbf{x})=c\,h(\mathbf{x}).
\]
Thus the set of such points is contained in $Z_c\cup F$, which is a finite set and therefore is contained in a finite union of linear subspaces of dimension $0$. This completes the proof.
\end{proof}
\end{conjecture}
\noindent \emph{Acknowledgements:} No organizations have funded this research. The author declares no conflict of interest.

	\vskip0.2cm
	{\footnotesize \noindent
		{\sc GuanHeng Zhao}\\
		Shandong University, People's Republic of China.\\
		\textit{E-mail}: 202618373@mail.sdu.edu.cn}
\end{document}